\documentclass[letterpaper]{template/mandc2019}
\usepackage{cancel}
\usepackage{tabls}
\usepackage{template/cites}
\usepackage{epsf}
\usepackage{xfrac}
\usepackage{appendix}
\usepackage{placeins}
\usepackage{pdflscape}
\usepackage{ragged2e}
\usepackage[top=1in, bottom=1.in, left=1.in, right=1.in]{geometry}
\usepackage{enumitem}
\setlist[itemize]{leftmargin=*}
\usepackage{caption}
\usepackage{subcaption}
\usepackage{algorithm}
\usepackage[noend]{algpseudocode}
\usepackage{amsmath,amssymb,amsfonts,amsthm}
\usepackage{color}
\usepackage[font=small,labelfont=bf]{caption}
\usepackage{bm}
\usepackage{lineno}

\usepackage[dvipsnames]{xcolor}

\usepackage{fancyvrb}

\RecustomVerbatimCommand{\VerbatimInput}{VerbatimInput}%
{fontsize=\footnotesize,
	framesep=2em, 
	rulecolor=\color{Gray},
	commandchars=\|\(\), 
	commentchar=*        
}

\newcommand{\SN}{S$_N$\,}
\newcommand{\PN}{P$_N$\,}

\newcommand{\aSN}{as-S$_N$\,}

\newcommand{\bOmega}{{\bm{\Omega}}}
\newcommand{\bx}{{\bm{x}}}

\title{Ray Effect Mitigation 
	for the Discrete Ordinates Method
	Using Artificial Scattering
}
\author{
	\textbf{
		Thomas Camminady${}^1$,
		Martin Frank${}^{1}$\footnote{Corresponding author. Submitted to Nuclear Science and Engineering on November 19, 2019.}$\,$, 
		Cory D. Hauck${}^2$, and 
		Jonas Kusch${}^1$
	}\\
	\texttt{thomas.camminady@kit.edu} \texttt{ martin.frank@kit.edu}\\
	\texttt{hauckc@ornl.gov} \texttt{ jonas.kusch@kit.edu}\\
	${}^1$ Karlsruhe Institute of Technology,\,
	${}^2$ Oak Ridge National Laboratory 
}
\begin{document}
\maketitle
\justify
\begin{abstract}
Solving the radiative transfer equation with the discrete ordinates (S$_N$) method leads to a non-physical imprint of the chosen quadrature set on the solution. 
To mitigate these so-called ray effects, we propose a modification of the \SN method, which we call artificial scattering \SN (as-S$_N$). 
The method adds an artificial forward-peaked scattering operator which generates angular diffusion to the solution and thereby mitigates ray effects. Similar to artificial viscosity for spatial discretizations, the additional term vanishes as the number of ordinates approaches infinity. 
Our method allows an efficient implementation of explicit and implicit time integration according to standard \SN solver technology. 
For two test cases, we demonstrate a significant reduction of the error for the as-S$_N$ method when compared to the standard S$_N$ method, both for explicit and implicit computations. 
Furthermore, we show that a prescribed numerical precision can be reached with less memory due to the reduction in the number of ordinates.
\end{abstract}

\section{Introduction} 

Several applications in the field of physics require an accurate solution of the radiative transfer equation. This equation describes the evolution of the angular flux of particles moving through a material medium. Examples include nuclear engineering \cite{duderstadt1976nuclear,henry1977nuclear}, high-energy astrophysics \cite{lowrie1999coupling,mcclarren2008manufactured}, supernovae \cite{fryer2006snsph,swesty2009numerical}, and fusion \cite{matzen2005pulsed,marinak2001three}. A major challenge when solving the radiative transfer equation numerically is the high-dimensional phase space on which it is defined. There are three spatial dimensions, two directional (angular) parameters, velocity, and time. In many applications, there is additional frequency or energy dependence. Hence, numerical methods to approximate the solution require a carefully chosen phase space discretization. 

There are several strategies to discretize the angular variables, and they all have certain strengths and weaknesses \cite{brunner2002forms}. The spherical harmonics (P$_N$) method 
\cite{case1967linear,pomraning1973equations,lewis1984computational} is a spectral Galerkin discretization of the radiative transfer equation. It uses the spherical harmonics basis functions to represent the solution in terms of 
angular variables with finitely many expansion coefficients, called moments. The P$_N$ method preserves rotational symmetry and shows spectral convergence for smooth solutions. However, like most spectral methods, P$_N$ yields oscillatory solution approximations in non-smooth regimes, which can lead to negative, non-physical angular flux values. Filtering of the expansion coefficients has shown to mitigate oscillations \cite{mcclarren2010robust} and a modified equation analysis has shown that filtering adds an artificial forward-peaked scattering operator to the equation if a certain scaling strength of the filter is chosen.


The discrete ordinates (S$_N$) method  
\cite{lewis1984computational} approximates the radiative transfer equation on a set of discrete angular directions.
The \SN discretization preserves positivity of the angular flux while yielding an efficient and straight forward implementation of time-implicit methods. 
However, the method is plagued by numerical artifacts, know as ray effects, when there are not enough ordinates to resolve the angular flux.
Because increasing the number of ordinates significantly increases numerical costs of simulations,
a major task to improve the solution accuracy of \SN methods is to mitigate these ray effects \cite{lathrop1968ray,morel2003analysis,mathews1999propagation} without simply adding more ordinates.

Various strategies to mitigate ray effects at affordable costs have been developed. 
In \cite{abu2001angular}, a biased quadrature set, which reflects the importance of certain ordinates is used. 
Furthermore, \cite{lathrop1971remedies} presents a method combining the \PN with the \SN method. Further studies for this method can be found in \cite{jung1972discrete,reed1972spherical,miller1977ray}, which show a reduction of ray effects. 
In \cite{morel2003analysis}, a comparison of these methods can be found.
In \cite{tencer2016ray}, computing the angular flux for differently oriented quadrature sets and averaging over different solutions has been proposed to reduce ray artifacts. 
In \cite{camminady2019ray}, a rotated \SN method has been developed, which rotates the quadrature set after every time iteration. Consequently, particles can move on a heavily increased set of directions of travel, leading to a reduction of ray effects. Analytic results show that rotating the quadrature set plays the role of an angular diffusion operator, which smears out artifacts that stem from the finite number of ordinates. Unfortunately, this method does not allow a straight forward implementation of sweeping, complicating the use of implicit methods.

The idea of this work is to add angular diffusion directly with the help of a forward-peaked artificial scattering operator. We choose this operator so that the effect of artificial scattering vanishes in the limit of infinitely many ordinates, but at finite order adds angular diffusion in such a way that it mitigates ray effects. Unlike the rotated \SN method in \cite{camminady2019ray}, the current approach allows for a straight forward implementation of sweeping, which we use to implement an implicit method.

\section{Main idea} 
In this section, we summarize the relevant mathematical background and introduce notation. We illustrate the problem of ray effects that occurs when discretizing the transport equation in angle and how artificial scattering can be used to mitigate these ray effects. We demonstrate that artificial scattering behaves like a Fokker-Planck operator in the appropriate limit.

\subsection{Radiative transfer equation}

The radiative transfer equation describes the evolution of the angular flux $\psi(t,\bx, \bOmega)$  via
\begin{linenomath*}\begin{equation}\label{eq:kineticEquation}
\partial_t \psi(t,\bx,\bOmega) + \bOmega \cdot \nabla_\bx \psi(t,\bx,\bOmega) + 
\sigma_{t}(\bx) \psi(t,\bx,\bOmega)= \sigma_{s}(\bx) (S^+ \psi) (t,\bx,\bOmega) + q(t,\bx),
\end{equation}\end{linenomath*}
where $t\in\mathbb{R}_+$ denotes time, $\bx\in\mathbb{R}^3$ is the spatial variable, and $\bOmega \in \mathbb{S}^2$ represents the direction.
The total cross section is $\sigma_{t}({\bx})=\sigma_{a}({\bx}) + \sigma_{s}({\bx})$. 
In the case of scattering, the in-scattering kernel operator $S^+(\psi)(t,\bx, \bOmega)$ describes the gain of particles that were previously traveling along direction $\bOmega'$ and changed to direction $\bOmega$. It is given by
\begin{linenomath*}\begin{equation}\label{eq:in-scattering}
(S^+ \psi) (t,\bx,\bOmega) = \int_{\mathbb{S}^2} s(\bOmega\cdot\bOmega')\psi(t,\bx,\bOmega') d\bOmega',
\end{equation}\end{linenomath*}
where $s(\bOmega\cdot\bOmega')$ is the probability of transitioning from direction $\bOmega'$ into direction $\bOmega$ or vice versa.
We assume---for simplicity---that the source $q(t,\bx)$ is isotropic.

\subsection{Ray effects}
\begin{figure}
	\centering
	\begin{subfigure}{.5\textwidth}
		\centering
	\includegraphics[width=0.99\linewidth]{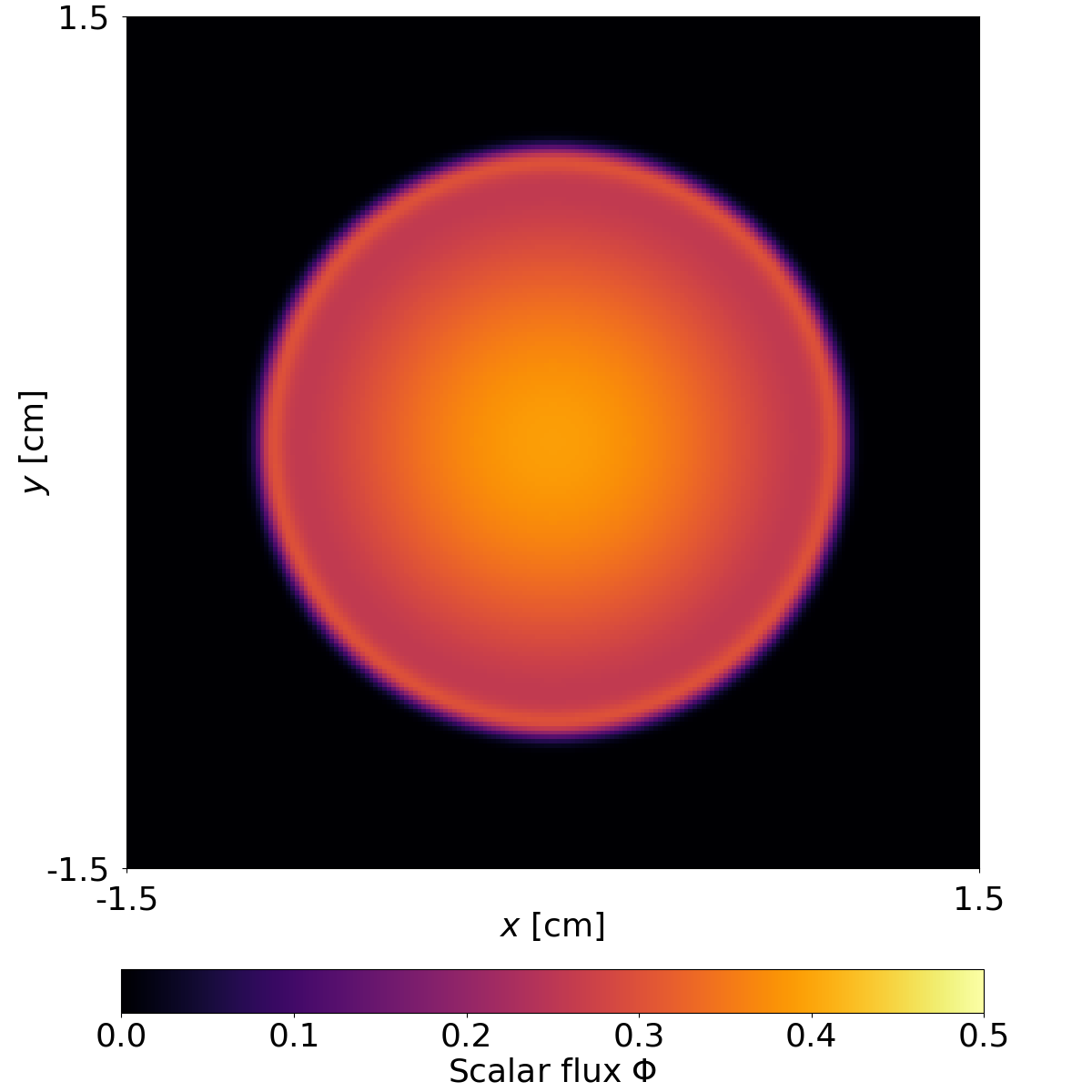}
		\caption{The semi-analytical reference solution to the line source problem.}
		\label{fig:semianalyticallinesource}
	\end{subfigure}%
	\begin{subfigure}{.5\textwidth}
		\centering
	\includegraphics[width=0.99\linewidth]{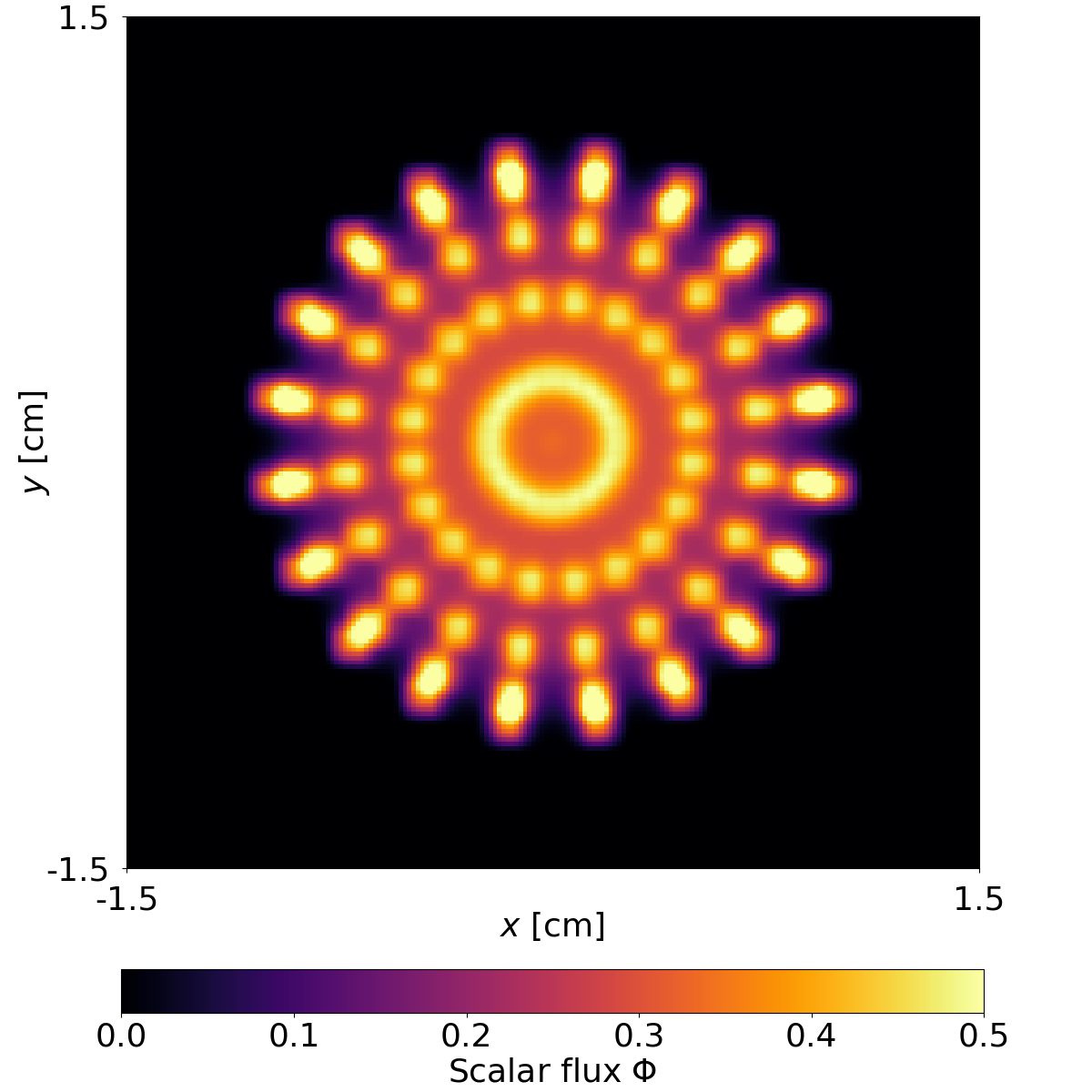}
		\caption{Numerical solution using the S$_N$ method with a tensorized quadrature. Ray effects are clearly visible.}
		\label{fig:numericallinesource}
	\end{subfigure}
	\caption{Illustrating ray effects with the line-source problem.}
	\label{fig:illustratingrayeffects}
\end{figure}

As previously explained,
the S$_N$ method preserves positivity but suffers from ray effects. 
An example of these artifacts is demonstrated for the line-source benchmark in Fig.\ref{fig:illustratingrayeffects}.
While the true scalar flux $\Phi(t,\bx):=\int_{\mathbb{S}^2}\psi(t,\bx,\bOmega') d\bOmega'$ is radially symmetric, the numerical solution has artifacts in the form of oscillations. We will discuss the line-source problem in more detail in Section \ref{subsec:results_linesource}.

\subsection{Artificial scattering}
We propose to address the problem of ray effects by adding artificial scattering to the right-hand side of \eqref{eq:kineticEquation}, in the form of an anisotropic scattering operator.
The modified system is then
\begin{linenomath*}\begin{equation}\label{eq:assneq}
\begin{split}
\partial_t \psi(t,\bx,\bOmega) &+ \bOmega \cdot \nabla_\bx \psi(t,\bx,\bOmega) + 
\sigma_{t}(\bx) \psi(t,\bx,\bOmega) \\
&= \sigma_{s}(\bx) (S^+ \psi) (t,\bx,\bOmega)+\sigma_{\text{as}}(\bx) (S_{\text{as}}\psi)(t,\bx,\bOmega) + q(t,\bx),
\end{split}
\end{equation}\end{linenomath*}
where
\begin{linenomath*}\begin{equation} \label{eq:artificialscattering}
	\sigma_{\text{as}}(\bx) (S_{\text{as}}\psi)(t,\bx,\bOmega) = \sigma_{\text{as}}(\bx) \int_{\mathbb{S}^2} s_\varepsilon(\bOmega' \cdot \bOmega) 
	\left(\psi(t,\bx,\bOmega')-\psi(t,\bx,\bOmega)\right)\, d\bOmega'.
	\end{equation}\end{linenomath*}

Here, $s_\varepsilon$ can be any Dirac-like sequence\footnote{While the idea of artificial scattering works with any Dirac-sequence, the asymptotic analysis that is performed later imposes stronger requirements to obtain a Fokker-Planck operator in the respective limit.} , i.e.
\begin{linenomath*}
\begin{align}
\int_{-1}^{+1} s_\epsilon(\mu)\, d\mu =1
 \text{ and }\int_{-1}^{+1}s_\varepsilon(\mu) f(\mu)\, d\mu \to f(1)
\end{align}
\end{linenomath*}
 for any sufficiently smooth function $f$ as $\varepsilon\to 0$. In our experiments, we choose
\begin{linenomath*}\begin{equation} \label{eq:scatteringkernel_artificialscattering}
s_\varepsilon(\mu) = \frac{2}{\sqrt{\pi}\, \varepsilon \,\text{Erf}\left(\frac{2}{\varepsilon}\right)}\, e^{-\sfrac{(1-\mu)^2}{\varepsilon^2}},
\end{equation}\end{linenomath*}
where the error function satisfies $\text{Erf}(x)\to 1$ as $x\to \infty$.
The proposed method, which we call artificial scattering-S$_N$ (as-S$_N$), has the following effects:
\begin{enumerate}
\item Similar to the artificial viscosity used to stabilize spatial discretizations of hyperbolic operators \cite[Chapter~16.1]{leveque1992numerical}, the artificial scattering adds an angular diffusion term to the radiative transfer equations. 
This term should vanish when the discretization is refined. Therefore, the variance of the artificial scattering kernel should be chosen to vanish in the limit $N_q\rightarrow\infty$.
We choose the variance to be $\varepsilon = \beta/N_q$, where $\beta$ is a constant, user-determined parameter. 
This choice ensures that $\varepsilon$ scales the average width of quadrature points, meaning that the domain of influence includes roughly the same number of ordinates when $N_q$ increases. In the limit, the \aSN solution converges to the classical \SN solution. 
\item The total number of particles is preserved by the artificial scattering term. Higher order moments are however dampened by the magnitude of artificial scattering. A beam of particles inside a void will be subject to scattering by the \aSN method, however artifacts that result from the standard \SN method dominate the overall error, unless the beam is aligned with the quadrature set. 
\item \aSN has similarities to the P$_{N-1}$-equivalent \SN method \cite{miller1977ray}: To mitigate ray effect, this method adds a fictitious source to the radiative transfer equation. This source, though derived by a different strategy, requires similar modifications of the standard \SN implementation. The main difference is that the artificial scattering kernel of \aSN is forward peaked, which can be used to design an efficient numerical treatment.
\item as-S$_N$ can be compared to filtered P$_N$ \cite{mcclarren2010robust,hauck2019filtered}, since the artificial scattering acts as a filter on the moment level.
\item The as-S$_N$ equation \eqref{eq:assneq} can---with appropriate boundary and initial conditions---be solved in a straight-forward manner using common S$_N$ implementations. When discussing one such implementation, we will focus on implicit discretization techniques and derive an efficient algorithm to treat the artificial scattering term.
\end{enumerate}

\subsection{Artificial scattering kernel}

To better distinguish between the two types of scattering, we will call the naturally occurring scattering of \eqref{eq:kineticEquation} \textit{physical scattering} and the scattering in \eqref{eq:artificialscattering} \textit{artificial scattering}. The way the artificial scattering is written in \eqref{eq:artificialscattering}, it includes in-scattering and out-scattering. We can split this further into 
\begin{linenomath*}\begin{equation} \label{eq:split_in_out_scattering}
(S_{\text{as}}\psi)(t,\bx,\bOmega) = \int_{\mathbb{S}^2} s_\varepsilon(\bOmega' \cdot \bOmega) 
\left(\psi(t,\bx,\bOmega')-\psi(t,\bx,\bOmega)\right)\, d\bOmega' = (S_{\text{as}}^+\psi)(t,\bx,\bOmega) -(S_{\text{as}}^-\psi)(t,\bx,\bOmega),
\end{equation}\end{linenomath*}
with 
\begin{linenomath*}\begin{equation}\label{eq:artificial_inscattering}
(S_{\text{as}}^+\psi)(t,\bx,\bOmega) = \int_{\mathbb{S}^2} s_\varepsilon(\bOmega' \cdot \bOmega) 
\psi(t,\bx,\bOmega')\, d\bOmega'
\end{equation}\end{linenomath*}
and
\begin{linenomath*}\begin{equation}\label{eq:artificial_outscattering}
(S_{\text{as}}^-\psi)(t,\bx,\bOmega) = \psi(t,\bx,\bOmega).
\end{equation}\end{linenomath*}

\subsection{Modified equation analysis}
According to Pomraning \cite{pomraning1992fokker}, the Fokker-Planck operator can be a legitimate description of highly peaked scattering.
This is true if (i) the scattering kernel $s_\epsilon(\mu)$ is a Dirac sequence, and (ii) the transport coefficients $p_{\varepsilon,i}:= \int_{-1}^{+1} (1-\mu)^i s_\varepsilon(\mu)\, d\mu$ are of order $\mathcal{O}(\varepsilon^i)$.
The resulting modified equation then reads
\begin{linenomath*}\begin{equation} \label{eq:fokkerplancklimit}
	\begin{split}
	\partial_t \psi(t,\bx,\bOmega) &+ \bOmega\cdot \nabla_\bx \psi + (\sigma_a+ \sigma_s ) \psi \\
	&= \sigma_s \cdot \Phi +\pi\cdot p_{\varepsilon,1}\cdot \sigma_{as}\cdot  \Delta_\bOmega\, \psi+ \mathcal{O}\left(\varepsilon^2\right),
	\end{split}
	\end{equation}\end{linenomath*}
where $\Delta_\bOmega$ is the Laplace operator in spherical coordinates.
We have already shown (i). 
To verify (ii), let $y=(1-\mu) / \varepsilon$. Then
\begin{linenomath*}\begin{align}
	p_{\varepsilon,i} &= \int_{2/\varepsilon}^0 (\varepsilon \, y)^i \frac{2}{\sqrt{\pi}\, \varepsilon \,\text{Erf}\left(\frac{2}{\varepsilon}\right)}
	e^{-y^2} (-\varepsilon) \, dy \\
	&=
	\frac{2}{\sqrt{\pi}\, \varepsilon \,\text{Erf}\left(\frac{2}{\varepsilon}\right)}\, \varepsilon^i \, \int_0^{2/\varepsilon}y^i e^{-y^2}\, dy \\
	&= \frac{2}{\sqrt{\pi}\, \varepsilon \,\text{Erf}\left(\frac{2}{\varepsilon}\right)}\, \varepsilon^i
	\left[\Gamma\left(\frac{1+i}{2}\right)-\Gamma\left(\frac{1+i}{2},\frac{4}{\varepsilon^2}\right)\right] \\
	&= \mathcal{O}(\varepsilon^i),
	\end{align}\end{linenomath*}
where $\Gamma(\cdot)$ and $\Gamma(\cdot,\cdot)$ denote the gamma function and the upper incomplete gamma function, respectively. 
Since (ii) implies that $p_{\varepsilon,1} = \mathcal{O}(\varepsilon)$, the operator vanishes if we let $\varepsilon\to 0$. 

We set $\varepsilon=\beta /N_q$ in the discrete case so that the angular diffusion vanishes if the number of ordinates $N_q$ tends to infinity. 
This analysis shows that the product $\sigma_{\text{as}}\cdot \beta$ controls the strength of the added angular diffusion. Section \ref{subsec:results_linesource} confirms this numerically.
\section{Discretization}
\label{sec:impl}
In the following section, we will discuss discretization and implementation of the presented \aSN method, laying the focus on how to incorporate artificial scattering into existing \SN codes.

\subsection{S$_N$ discretization}
For sake of completeness, we briefly summarize the \SN method. Given a finite number of ordinates $\bOmega_1,\dots,\bOmega_{N_q}$ and defining the \SN solution $\psi_q(t,\bx) \approx \psi(t,\bx,\bOmega_q)$ the \SN method solves the semi-discretized system of $N_q$ equations
\begin{linenomath*}\begin{equation}
\begin{split}
\label{eq:SN}
\partial_t \psi_q(t,\bx) + \bOmega_q \cdot \nabla_\bx \psi_q(t,\bx) + \sigma_{t}(\bx) \psi_q(t,\bx)= \sigma_{s}(\bx) \sum_{p=1}^{N_q} w_p \cdot s(\bOmega_q\cdot\bOmega_p) \psi_p(t,\bx) + q(t,\bx).
\end{split}
\end{equation}\end{linenomath*}
Here, $w_p$ are quadrature weights, chosen such that 
\begin{linenomath*}\begin{equation}
\int_{\mathbb{S}^2} \psi(t,{\bx}, {\bOmega}) \approx \sum_{q=1}^{N_q} w_q \cdot \psi(t,\bx,\bOmega_q).
\end{equation}\end{linenomath*}
To compute numerical solutions, we still need to discretize \eqref{eq:SN} in space and time. In this work, we will investigate solutions for both implicit and explicit time discretizations. The explicit code uses Heun's method as well as a minmod slope limiter. It is based on \cite{camminady2019ray,garrett2013comparison}, which provide a detailed description of the chosen methods. The implicit discretization and an efficient strategy to integrate artificial scattering in a given implicit code framework will be discussed in Sections~\ref{sec:Implicit} and \ref{sec:ImplicitImplementation}.

\subsection{Adding artificial scattering to the \SN equations}
\label{sec:implassn}
Our goal is to include artificial scattering in the S$_N$ equations in \eqref{eq:SN}. By simply approximating the artificial scattering term in \eqref{eq:assneq} with the chosen quadrature rule, we obtain the \aSN equations 
\begin{linenomath*}\begin{equation}
\begin{split}
	\label{eq:asSNeq}
    \partial_t \psi_q(t,\bx) + \bOmega_q &\cdot \nabla_\bx \psi_q(t,\bx) + \sigma_t(\bx) \psi_q(t,\bx) + \sigma_{\text{as}}(\bx)\psi_q(t,\bx)
    \\=&\sigma_s(\bx) \sum_{p=1}^{N_q} w_p \cdot c_q \cdot s(\bOmega_q\cdot\bOmega_p) \psi_p(t,\bx)
    \\
    &+\sigma_{\text{as}}(\bx)\sum_{p=1}^{N_q}w_p \cdot c_q^{(\varepsilon)}\cdot s_\varepsilon(\bOmega_q\cdot \bOmega_p)\psi_p(t,\bx) \\
    &+  q(t,\bx) 
    .
    \end{split}
\end{equation}\end{linenomath*}
Here, $c_q := 1 / \sum_p w_p \cdot s(\bOmega_q\cdot \bOmega_p)$ and $c_q^{(\varepsilon)} := 1 / \sum_p w_p \cdot s_\varepsilon(\bOmega_q\cdot \bOmega_p)$ are normalization factors. While on the continuous level, these factors are the same for every direction, we obtain a dependency on the chosen ordinate due to the non-uniform discretization in angle. These normalization factors are needed to obtain a simple expression for the out-scattering terms. Moving these terms to the left-hand side of \eqref{eq:asSNeq} stabilizes the source iteration used in the implicit method.

\subsection{Quadrature}
\label{subsec:quadrature}
\begin{figure}[h!]
	\centering
	\includegraphics[scale=0.4]{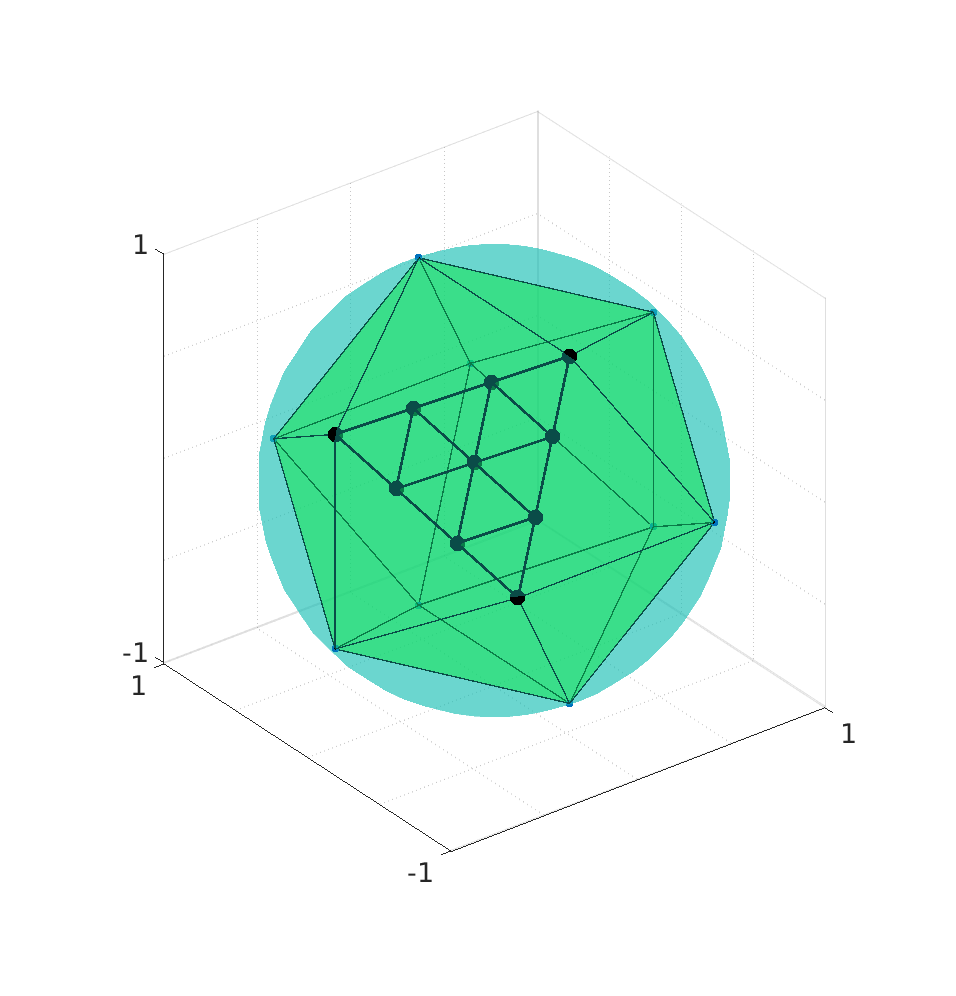}
	\caption{Construction of the icosahedron quadrature set. One can see the icosahedron geometry with a triangulated face which is then mapped onto the sphere. For the quadrature set, we align one of the vertices with the point $(0,0,1)$.}
	\label{fig:icosahedronQuadrature}
\end{figure}

It remains to pick an adequate quadrature set. When applying artificial scattering, the solution smears out along the imprint of the discrete set of ordinates. To ensure an evenly spread artificial scattering effect, a quadrature with a highly uniform ordinate distribution should be chosen. Commonly, the construction of a quadrature set starts at a chosen planar geometry, which is discretized and then mapped onto the surface of a sphere. The mapped nodes of the previously chosen discretization are then taken to be the quadrature points while the weights are determined by the area associated with these points. An even node distribution is achieved by using an icosahedron as initial planar geometry, which one can find in Figure~\ref{fig:icosahedronQuadrature}. Each face of the icosahedron is triangulated to generate the nodes which will be mapped onto the surface of the sphere. There are different strategies to perform this triangulation and we choose an equidistant spacing of points on each line of the triangle as well as the corresponding points inside the triangle. The corresponding weight is the area of the hexagon which lies around the given node and is defined by connecting the midpoints of the neighboring triangles. For more details on the icosahedron quadrature, see \cite{icosahedron}.
From now on, we will exclusively use this quadrature in all S$_N$ and as-S$_N$ computations.

\subsection{Implicit time discretization}
\label{sec:Implicit}
\begin{figure}
	\centering
	\includegraphics[width=0.99\linewidth]{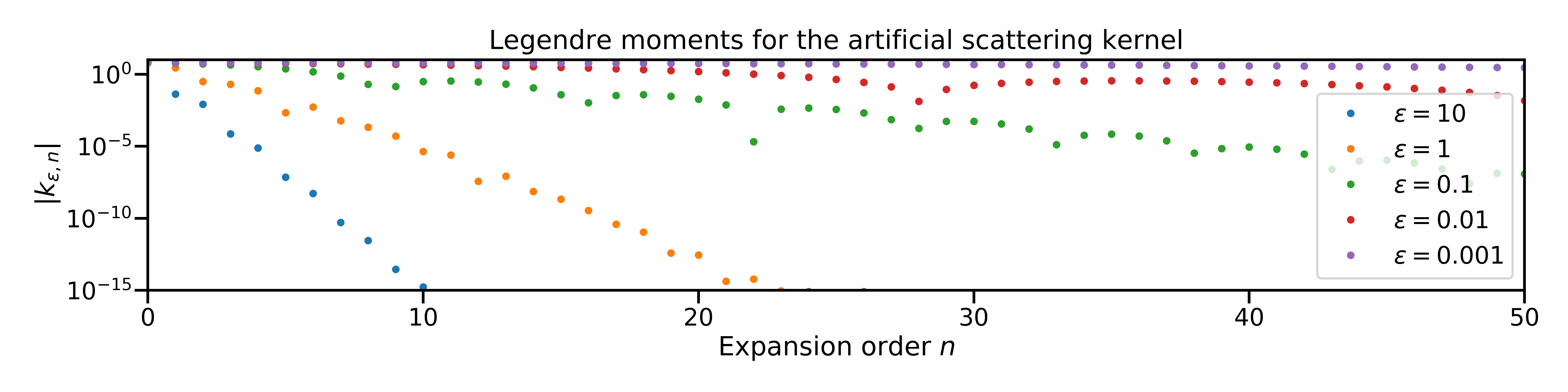}
	\caption{Decay of the Legendre moments $k_{\varepsilon,n} = 2\pi \int_{-1}^{+1} s_\varepsilon(\mu)P_n(\mu)\, d\mu$ for different values of $\varepsilon$ and the expansion order $n$.}
	\label{fig:kerneldecay}
\end{figure}

Implicit time discretization methods provide stability for large time steps, which are crucial in applications involving different time scales. However, when discretizing the radiative transfer equation, they require a matrix solve in every time step, which is commonly performed by a Krylov solver \cite{faber1988look,ashby1991preconditioned}. We start with an implicit Euler discretization, where we, in an abuse of notation, denote the flux at the new time step by $\psi(\bx,\bOmega)$ and at the old time step by $\psi^{\text{old}}(\bx,\bOmega)$. The equivalent \aSN system is
\begin{linenomath*}
\begin{align}
\bOmega \cdot \nabla_\bx \psi + \left(\sigma_a+ \sigma_s + \sigma_{\text{as}} + \frac{1}{\Delta t} \right)\psi = \sigma_s S^+\psi+\sigma_{\text{as}}S_{\text{as}}^+\psi+q + \frac{\psi^{\text{old}}}{\Delta t}.
\end{align}
\end{linenomath*}
Defining the streaming operator $L\psi := \bOmega \cdot \nabla_\bx \psi + \left(\sigma_a+ \sigma_s + \sigma_{\text{as}} + \frac{1}{\Delta t} \right)\psi$ as well as the modified source $\tilde q := q + \psi^{\text{old}} / \Delta t$, we can put this into more compact notation
\begin{linenomath*}\begin{align}
\label{eq:asSNImplCompact}
L\psi = \sigma_s S^+\psi+\sigma_{\text{as}}S_{\text{as}}^+\psi+\tilde q.
\end{align}\end{linenomath*}
First, let us numerically treat the artificial scattering in the same way as commonly done for physical scattering. The physical in-scattering kernel can be written as
\begin{align*}
S^+ = O \Sigma M,
\end{align*}
where $\Sigma$ carries the respective expansion coefficients of the scattering kernel, $M$ maps from the ordinates to the moments and $O$ from the moments back to the angular space. Making use of this strategy to represent the artificial scattering, we get
\begin{linenomath*}\begin{align}
S_{\text{as}}^+ = O\Sigma_{\text{as}}M.
\end{align}\end{linenomath*}
When denoting the moments as $\phi = M\psi$, equation \eqref{eq:asSNImplCompact} becomes
\begin{linenomath*}
\begin{align}
L\psi = \sigma_s O \Sigma\phi +\sigma_{\text{as}}O \Sigma_{\text{as}}\phi+\tilde q\;.
\end{align}
\end{linenomath*}
Inverting $L$ and applying $M$ to both sides yields the fixed point equations
\begin{linenomath*}
	\begin{align}\label{eq:KrylovNaive}
\phi = \sigma_s M L^{-1} O \Sigma\phi +\sigma_{\text{as}}ML^{-1}O\Sigma_{\text{as}} \phi+ML^{-1}\tilde q \;.
\end{align}
\end{linenomath*}
Note that with $\sigma_{\text{as}}=0$, this is the standard equation to which a Krylov solver is applied. Choosing a non-zero artificial scattering strength can result in significantly increased numerical costs when solving \eqref{eq:KrylovNaive} with a Krylov method: To show this, let us move to the discrete level, i.e. discretizing the directional domain, which requires picking a finite number of moments. In this case $\Sigma$ becomes a diagonal matrix with entries falling rapidly to zero (in the case of isotropic scattering, only the first entry is non-zero). Hence, few moments are required to capture the effects of physical scattering. However, since the artificial scattering kernel is strongly forward-peaked, the entries of the diagonal matrix $\Sigma_{\text{as}}$ do not fall to zero quickly, meaning that the method requires a large number of moments to include artificial scattering, which results in a heavily increased run time \cite{hauck2019filtered}. The slow decay of the Legendre moments $k_{\varepsilon,n} = 2\pi \int_{-1}^{+1} s_\varepsilon(\mu)P_n(\mu)\, d\mu$ for $\varepsilon\rightarrow 0$ is visualized in Fig.~\ref{fig:kerneldecay}.
In order to be able to choose the reduced number of moments required to resolve physical scattering, we move the artificial scattering into the sweeping step. Hence, going back to equation \eqref{eq:asSNImplCompact}, we only perform the moment decomposition on the physical scattering to obtain
\begin{linenomath*}\begin{align}
(L - \sigma_{\text{as}}S_{\text{as}}^+)\psi = \sigma_s O \Sigma\phi + \tilde q\;.
\end{align}\end{linenomath*}
Moving the operator $(L - \sigma_{\text{as}}S_{\text{as}}^+)$ to the right hand side and taking moments yields
\begin{linenomath*}\begin{align}\label{eq:KrylovFixedPoint}
\phi = \sigma_s M(L - \sigma_{\text{as}}S_{\text{as}}^+)^{-1} O \Sigma\phi + M(L - \sigma_{\text{as}}S_{\text{as}}^+)^{-1}\tilde q\;.
\end{align}\end{linenomath*}
The Krylov solver is then applied to this fixed-point iteration. In contrast to \eqref{eq:KrylovNaive}, the physical scattering term dictates the number of moments. 
The computation of $(L-\sigma_{\text{as}}S_{\text{as}}^+)^{-1}$ is performed by a source iteration, where the general equation $(L-\sigma_{\text{as}}S_{\text{as}}^+)\psi = R$ is solved by
iterating on
\begin{linenomath*}\begin{align}\label{eq:sourceIteration}
L\psi^{(l+1)} = \sigma_{\text{as}}S_{\text{as}}^+\psi^{(l)} + R.
\end{align}\end{linenomath*}
This iteration is expected to converge fast since effects of artificial scattering will be small in comparison to physical scattering. 

\subsection{Implementation details}
\label{sec:ImplicitImplementation}
At this point, we choose a finite number of ordinates and moments, i.e. the flux $\psi$ is now a vector with dimension $N_q$ and the moments $\phi$ have finite dimension $N$. Consequently, operators applied to the directional space become matrices. For better readability, we abuse notation and reuse the same symbols as before. 

We observed that a second-order spatial scheme is required to capture the behavior of the test cases used in this work. To ensure an efficient sweeping step, we use a second-order upwind stencil without a limiter. Let us denote the operator $L$ discretized in space and direction by $L_{\Delta}$. For ease of presentation, we assume a slab geometry. i.e. we have the spatial variable $x\in\mathbb{R}$ and the directional variable $\mu\in[-1,1]$. In the following, we split the directional variable into $\mu_-\in[-1,0]$ and $\mu_+\in(0,1]$. An extension to arbitrary dimension is straight forward. Now with $\lambda_{\pm}:= \mu_{\pm}\frac{ \Delta t}{\Delta x}$ and $\sigma_t := \sigma_a+ \sigma_s + \sigma_{\text{as}} + \frac{1}{\Delta t}$, we can write the discretized streaming operator as
\begin{linenomath*}\begin{align}
L_{\Delta}\psi := \lambda_{\pm} ( g_{j+1/2} - g_{j-1/2} ) + \Delta t \sigma_t \psi.
\end{align}\end{linenomath*}
The numerical flux for $\mu_+$ is then given by
\begin{linenomath*}\begin{align}
g_{j+1/2} := a \psi_j + b \psi_{j-1}, \qquad \text{ with } a := \frac32, \enskip b:= -\frac12
\end{align}\end{linenomath*}
and for $\mu_-$ by
\begin{linenomath*}\begin{align}
g_{j+1/2} := a \psi_{j+1} + b \psi_{j+2}.
\end{align}\end{linenomath*}
This scheme is L$^2$ stable, which we show in Appendix~\ref{app:upwind}. Let us now discuss the implementation of the implicit method in more detail. As mentioned earlier, a source iteration \eqref{eq:sourceIteration} is required to invert the operator $(L-\sigma_{\text{as}}S_{\text{as}}^+)$. For an initial guess $\psi^{(0)}$ and an arbitrary right hand side $R$, this iteration is given by Alg.~\ref{alg:SourceIteration}. Note that the discrete artificial in-scattering $S_{\text{as}}^+$ is a sparse matrix, which guarantees an efficient evaluation of the matrix vector product $S_{\text{as}}^+\psi^{(l)}$ in \eqref{eq:sourceIteration}. Furthermore, the inverse of $L_{\Delta}$ can be computed by a sweeping procedure.
\begin{algorithm}[H]
\begin{algorithmic}[1]
\Procedure{SourceIteration}{$\psi^{(0)},R$}
\State $\ell \leftarrow 0$ 
\State $\psi^{(\ell+1)} \leftarrow L_{\Delta}^{-1}\left(\sigma_{\text{as}} S_{\text{as}}^+ \psi^{(\ell)} + R\right)$
  \While{$\Vert\psi^{(\ell+1)} - \psi^{(\ell)}\Vert_2\geq\epsilon\tilde c$}
\State $\psi^{(\ell+1)} \leftarrow L_{\Delta}^{-1}\left(\sigma_{\text{as}} S_{\text{as}}^+ \psi^{(\ell)} + R\right)$
\State $\ell \leftarrow \ell+1$
\EndWhile
\Return $\psi^{(\ell)}$
\EndProcedure
\end{algorithmic}
\caption{Source Iteration algorithm}
\label{alg:SourceIteration}
\end{algorithm}
In order to get a good error estimator, we set the constant $\tilde c := (1-T)/T$ in Algorithm~\ref{alg:SourceIteration}, where $T$ is an estimate of the Lipschitz constant and $\epsilon$ is a user-determined parameter. Our implementation solves the linear system of equations
\begin{linenomath*}\begin{subequations} 
\begin{align}
 A \phi &= b \label{eq:KrylovFixedPoint1} \\
\text{with }A &:=  - \sigma_s M(L - \sigma_{\text{as}}S_{\text{as}}^+)^{-1} O \Sigma \\
b &:= M(L - \sigma_{\text{as}}S_{\text{as}}^+)^{-1}\tilde q
\end{align}
\end{subequations} \end{linenomath*}
using a GMRES solver. The solver requires the evaluation of the left-hand side for a given $\psi$ with an initial guess $\psi^{(0)}$, which is given by Alg.~\ref{alg:LinearFunction}.
\begin{algorithm}[H]
\begin{algorithmic}[1]
\Procedure{LHS}{$\psi^{(0)},\phi$}
\State $\tilde\psi \leftarrow \text{SourceIteration}(\psi^{(0)},\sigma_s O\Sigma \phi)$
\State \Return $\phi - M\tilde\psi$
\EndProcedure
\end{algorithmic}
\caption{Left-hand side of \eqref{eq:KrylovFixedPoint1}}
\label{alg:LinearFunction}
\end{algorithm}
The main time stepping scheme is then given by Alg.~\ref{alg:Sweeping-Krylov}. After initializing $\psi$ and $\phi$, the right hand side to \eqref{eq:KrylovFixedPoint1} is set up in line 4. Line 5 then solves the linear system \eqref{eq:KrylovFixedPoint1} and Line 6 determines the time-updated flux $\psi$ from the moments $\phi^{\text{new}}$.
\begin{algorithm}[H]
\begin{algorithmic}[1]
\State $\psi^{\text{old}}\leftarrow \text{InitialCondition}()$
\State $\phi^{\text{old}} \leftarrow M\psi^{\text{old}}$
\While{$t<t_{\text{end}}$}
\State $b \leftarrow M \cdot \text{SourceIteration}(\psi^{old},q + \frac{1}{\Delta t}\psi^{old})$
\State $\phi^{\text{new}} \leftarrow \text{Krylov}(\text{LHS}(\psi^{\text{old}},\phi^{\text{old}}),b)$
\State $\psi^{\text{new}} \leftarrow \text{SourceIteration}(\psi^{\text{old}},\sigma_s O \Sigma \phi^{\text{new}} + \frac{1}{\Delta t}\psi^{\text{old}})$
\State $\psi^{\text{old}} \leftarrow \psi^{\text{new}}$
\State $\phi^{\text{old}} \leftarrow \phi^{\text{new}}$
\EndWhile
\end{algorithmic}
\caption{Sweeping-Krylov algorithm}
\label{alg:Sweeping-Krylov}
\end{algorithm}

There exist several ways to modify the presented algorithm to achieve higher performance. For example, 
one can modify the presented method by not fully converging the source iteration in Alg.~\ref{alg:SourceIteration}. Instead, only a single iteration can be performed to drive the moments $\phi$ and the respective angular flux $\psi$ to their corresponding fixed points simultaneously. In numerical tests, we observe that this will significantly speed up the calculation. However, since we do not focus on runtime optimization, we do not further discuss this idea and leave it to future work.

\section{Results}

In the following, we evaluate the proposed method within the scope of two numerical test cases: (i) the line-source problem is used as it is inherently prone to ray-effects when using the S$_N$ method, and (ii) the lattice test case models---in a very simplified way---neutrons in a fission reactor with a source and heterogeneous materials. For both problems, we present results for the explicit and implicit methods, respectively. 

Both test cases are computed on a two-dimensional regular grid for the spatial variable. We project the $\bOmega_q \in \mathbb{S}^2$ for $q=1,\ldots,N_q$ onto the $x$-$y$-plane. 

The code used to compute the numerical results is published under the MIT license in a public repository at \texttt{https://github.com/camminady/SN}.

\label{sec:res}
\subsection{Line-source test case}
\label{subsec:results_linesource}
The goal of this test case is to numerically compute the Green's function for an initial isotropic Dirac-mass at the origin, i.e.\ $\psi(t=0,\bx,\bOmega) = \sfrac{1}{4\pi}\,\delta(\bx)$, which is realized as a narrow Gaussian in space with $\psi(t=0,\bx,\bOmega) = \max\{10^{-4},\sfrac{1}{4\pi \delta}\,\exp( \sfrac{-\bx^2}{4\delta})\}$ and $\delta = 0.03^2$. We choose $\sigma_s = \sigma_t = 1$. The spatial discretization varies from $50\times 50$ for a coarse grid to $200\times 200$ points on the domain $[-1.5,1.5]\times [-1.5,1.5]$ for a fine grid. There exists a semi-analytical solution to the full transport equation for this problem due to Ganapol et al.~\cite{ganapol2001homogeneous}. The exact solution consists of a circular front moving away from the origin as well as a tail of particles which have been scattered or not emitted perpendicularly from the center. We chose the line-source because it is a test case that lays bare almost any artifact an angular discretization might suffer from.

The parameters of the artificial scattering have been set to 
$\varepsilon = \beta / N_q$ with $N_q$ as the number of quadrature points. 
Obviously, choosing these parameters requires some experience. However, as in the case of filtering for $P_N$ \cite{mcclarren2010robust}, both parameters can be adjusted for coarse angular and spatial grids, and are expected to be valid for finer grids, which seems to be the case for the line-source problem as well.
The CFL number, i.e. the ratio of the time step and the spatial cell size is $0.95$ for the explicit calculations and $2$ for the implicit calculation.
For the implicit discretization, the tolerance for the GMRES solver was set to $1.5\cdot10^{-8}$, and we considered the inner source iteration to be converged at an estimated error of $10^{-4}$.

\begin{figure}[h!]
	\centering
	\includegraphics[width=0.99\linewidth]{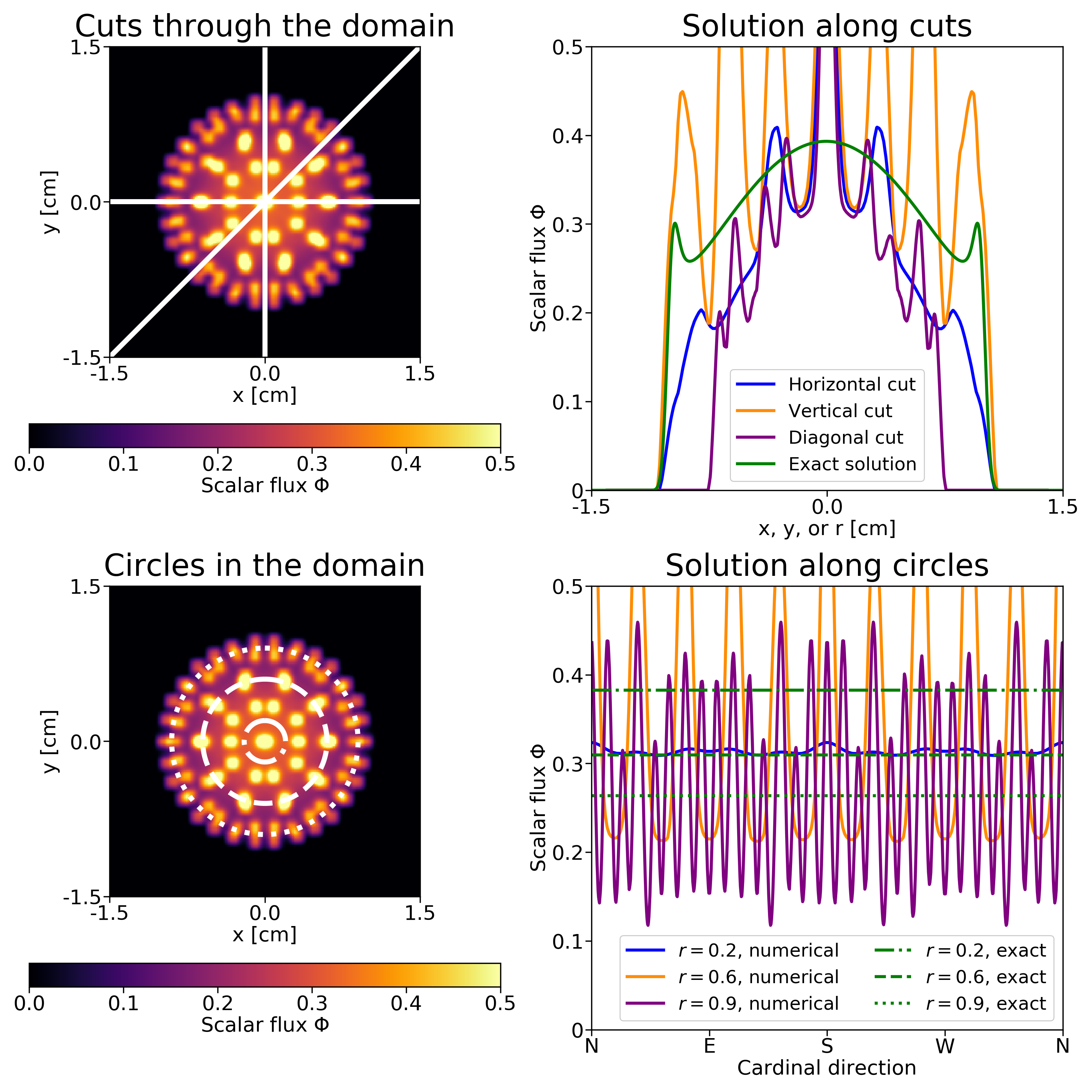}
	\caption{$S_4$ solution with ray-effects. We choose $N_q=92$ quadrature points, the spatial domain is composed of $N_x \times N_y = 200\times 200$ spatial cells and the CFL number is 0.95. Cuts through the domain and along circles with different radii are visualized in the right column. Only the solution along the horizontal cut is symmetric for the icosahedron quadrature.}
	\label{fig:1summarysn}
\end{figure}

An overview of the analytics that we perform is given in Fig.~\ref{fig:1summarysn}. We evaluate the scalar flux $\Phi(t,\bx) = \int_{\mathbb{S}^2} \psi(t,\bx,\bOmega')\,d\bOmega'$ at the final time step.
We have performed an explicit S$_4$ computation with $N_q=92$ and $N_x \times N_y = 200\times 200$. In both rows of Fig.~\ref{fig:1summarysn}, the left column shows the scalar flux at the final time step. The first row shows the solution along cuts through the domain on the right with the respective cuts on the left in white. The second row shows the solution along circles with different radii on the right and the respective circles on the left. Strong oscillations are visible due to ray-effects. For the first row, the analytical solution is given in green in the right column image. In the lower row, the analytical solution is constant along a circle with a certain radius, visualized for $r=0.2$, $r=0.6$, and $r=0.9$ in green.

\begin{figure}
	\centering
	\includegraphics[width=0.99\linewidth]{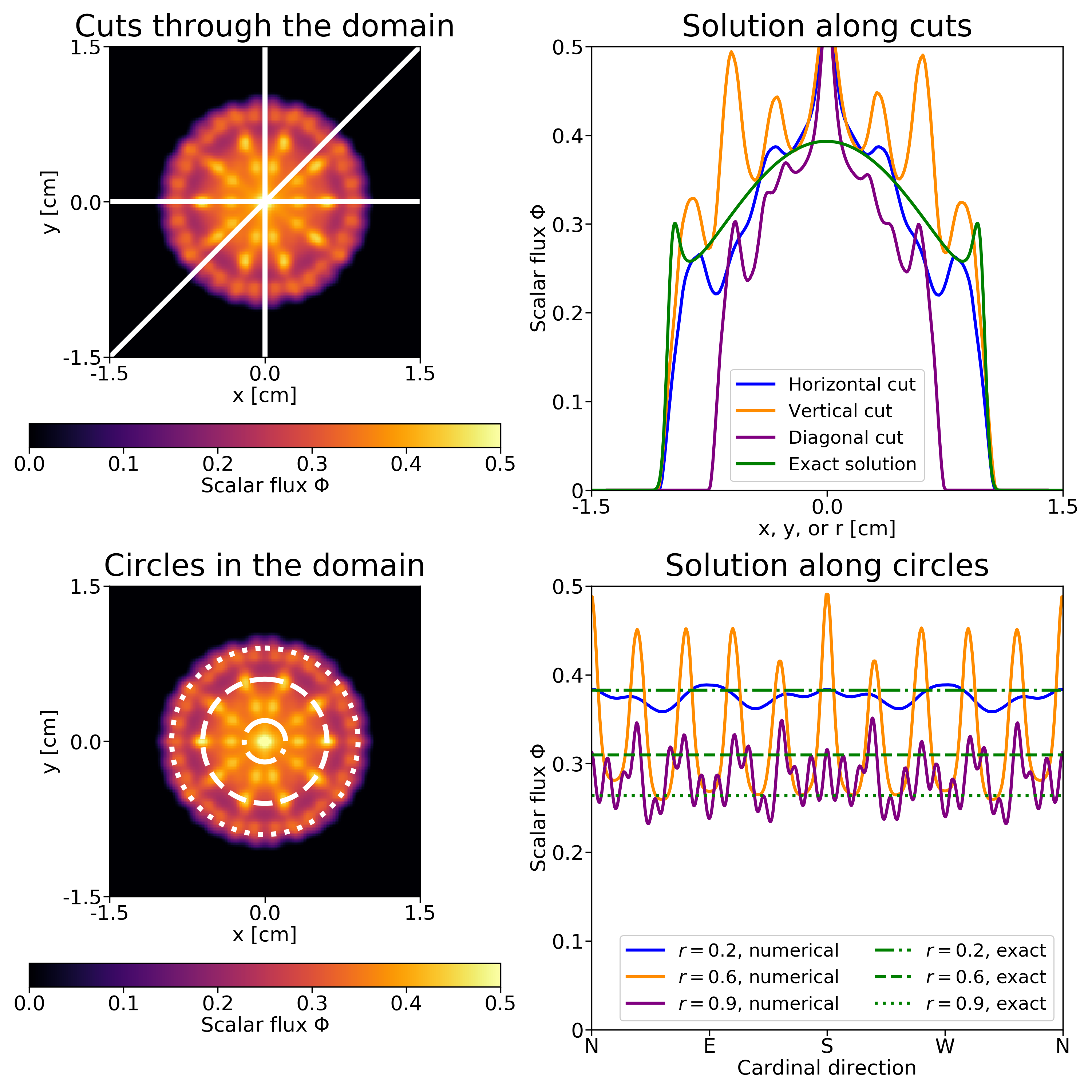}
	\caption{as-$S_4$ solution with mitigated ray-effects. We choose $N_q=92$ quadrature points, the spatial domain is composed of $N_x \times N_y = 200\times 200$ spatial cells and the CFL number is 0.95. Cuts through the domain and along circles with different radii are visualized in the right column.
		We set $\sigma_{\text{as}}=5$ and $\beta=4.5$. Only the solution along the horizontal cut is symmetric for the icosahedron quadrature.}
	\label{fig:1summaryassn}
\end{figure}

\begin{figure}
	\centering
	\includegraphics[width=0.99\linewidth]{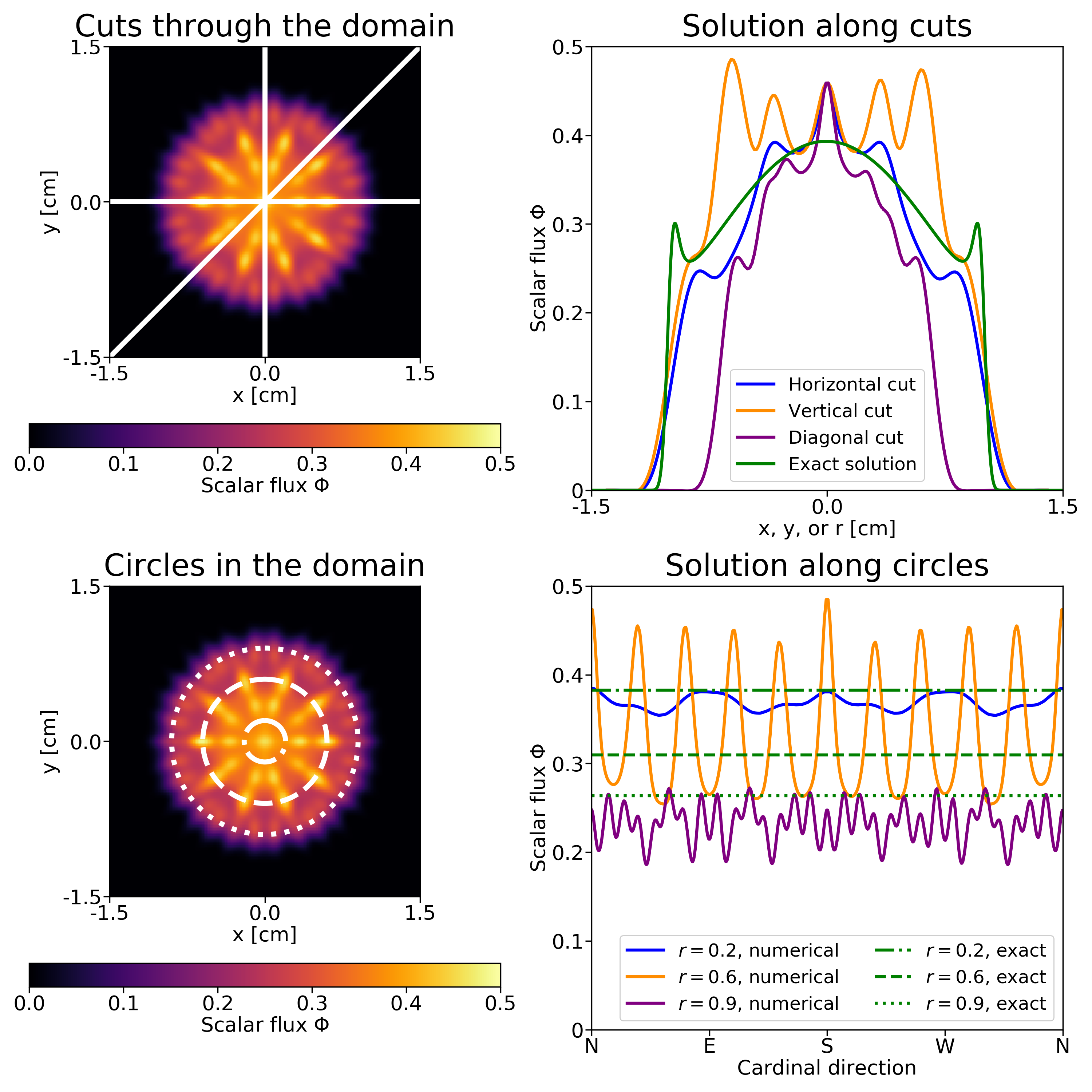}
			\caption{as-$S_4$ solution with mitigated ray-effects. We choose $N_q=92$ quadrature points, the spatial domain is composed of $N_x \times N_y = 200\times 200$ spatial cells and the CFL number is 2. Cuts through the domain and along circles with different radii are visualized in the right column. We set $\sigma_{\text{as}}=7$ and $\beta=4$. Only the solution along the horizontal cut is symmetric for the icosahedron quadrature.}
	\label{fig:1summaryassnimpl}
\end{figure}

In Fig.~\ref{fig:1summaryassn} and Fig.~\ref{fig:1summaryassnimpl} we see the same summary of results, now for an explicit and implicit computation, respectively. In both computations, ray-effects have been reduced significantly when comparing the results with Fig.~\ref{fig:1summarysn} despite the same number of quadrature points. The implicit calculation looks slightly more diffusive. However, the line-source problem is not a problem that would be computed implicitly in the first place and we only use it to illustrate the expected behavior for implicit computations.

\begin{landscape}
\begin{figure}
	\centering
	\includegraphics[width=0.99\linewidth]{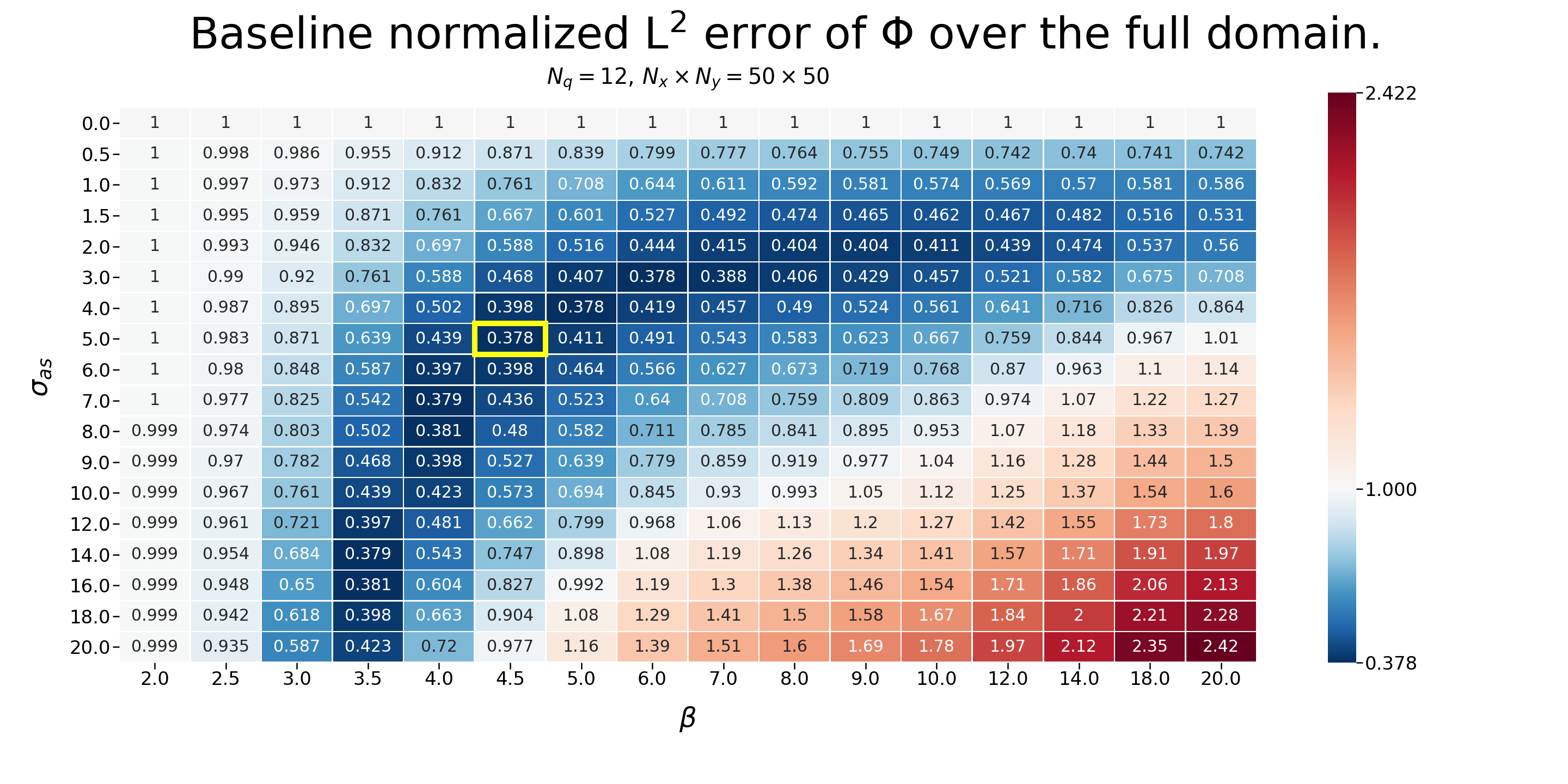}
	\caption{Parameter study for $\sigma_{\text{as}}$ and $\varepsilon=\beta/N_q$ on a grid of $N_q\times N_x \times N_y = 12\times50\times50$ in an explicit calculation.
	For every simulation we compute the L${}^2$ error of the scalar flux $\Phi$ with respect to a semi-analytical reference solution on the same spatial grid. The number in each field of the heatmap is then the baseline normalized error, i.e. the L${}^2$ error obtained for that specific parameter configuration divided by the error obtained without artificial scattering. For the case of $\beta=4.5$ and $\sigma_{\text{as}}=5$ (highlighted in yellow) the error drops down to $37.8\%$ of the original error without artificial scattering.}
	\label{fig:heatmapabsl1normalized}
\end{figure}
\end{landscape}

\begin{landscape}
	\begin{figure}
		\centering
	\includegraphics[width=0.99\linewidth]{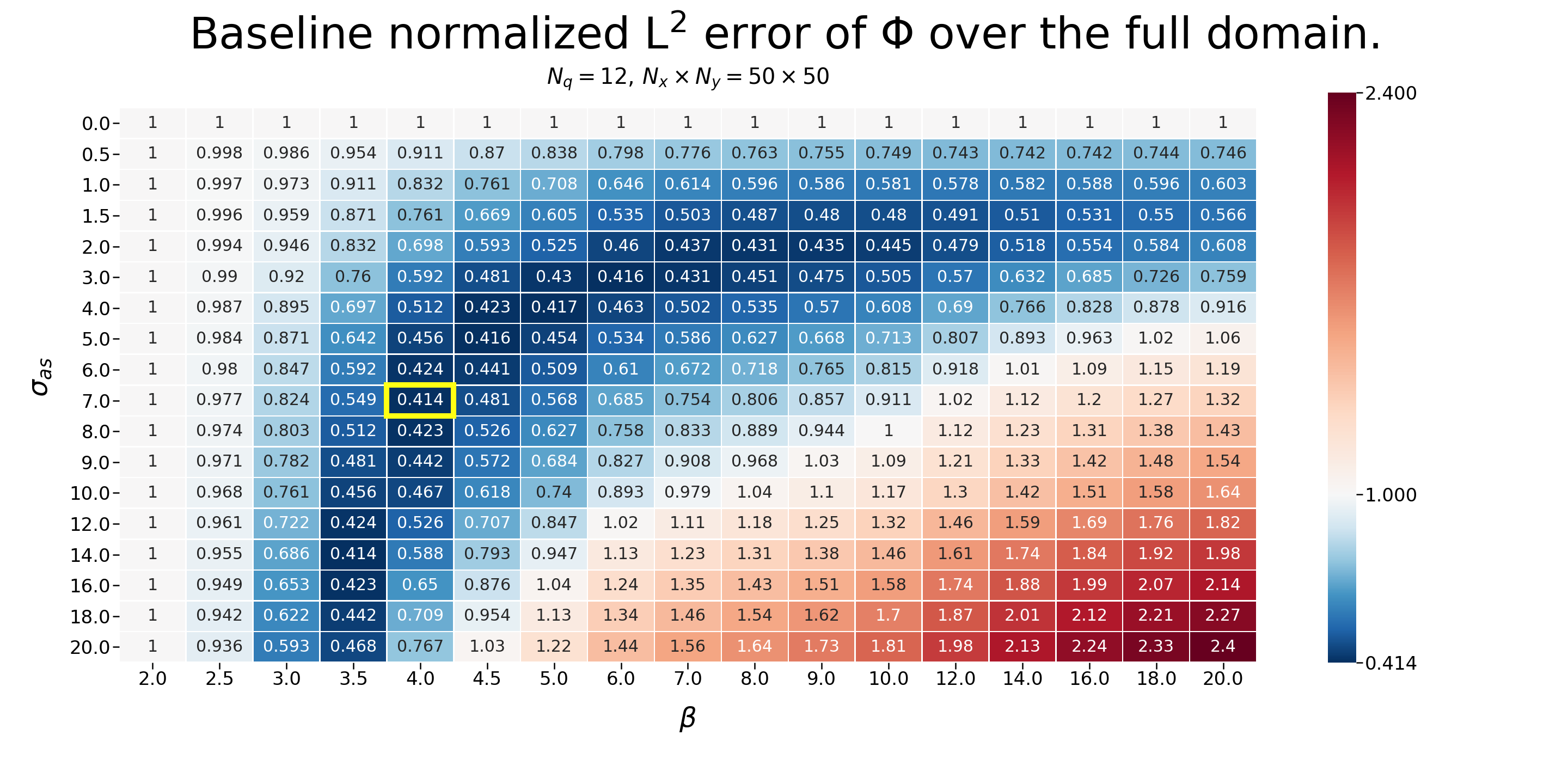}
	\caption{Parameter study for $\sigma_{\text{as}}$ and $\varepsilon=\beta/N_q$ on a grid of $N_q\times N_x \times N_y = 12\times50\times50$ in an implicit  calculation.
	For every simulation we compute the L${}^2$ error of the scalar flux $\Phi$ with respect to a semi-analytical reference solution on the same spatial grid. The number in each field of the heatmap is then the baseline normalized error, i.e. the L${}^2$ error obtained for that specific parameter configuration divided by the error obtained without artificial scattering. For the case of $\beta=4$ and $\sigma_{\text{as}}=7$ (highlighted in yellow) the error drops down to $41.4\%$ of the original error without artificial scattering.}%
		\label{fig:heatmapabsl1normalizedimpl}
	\end{figure}
\end{landscape}

The values for $\beta$ and $\varepsilon$ in the fine calculations are determined from a parameter study using coarse spatial and angular grids. The results of this parameter study are given in Fig.~\ref{fig:heatmapabsl1normalized} for the explicit algorithm and in Fig.~\ref{fig:heatmapabsl1normalizedimpl} for the implicit algorithm.


A single simulation for the coarse configuration takes $\sim \sfrac{1}{400}$ times the time of a single computation for the fine configuration. Consequently, the full parameter study with all $306$ configurations can be performed for less than the costs of a single fine computation. For the optimal parameter configuration, the error decreases down to $37.8\%$ for the explicit case and down to $41.4\%$ for the implicit case.

In both cases, implicit and explicit, we observe a region of parameters that yield similarly good results. This behavior mostly matches the predicted relation from the asymptotic analysis, i.e. when $\varepsilon$ is small, $\sigma_{\text{as}}\cdot \varepsilon$ controls the effect of artificial scattering.



\begin{figure}[h!]
	\centering
	\includegraphics[width=0.99\linewidth]{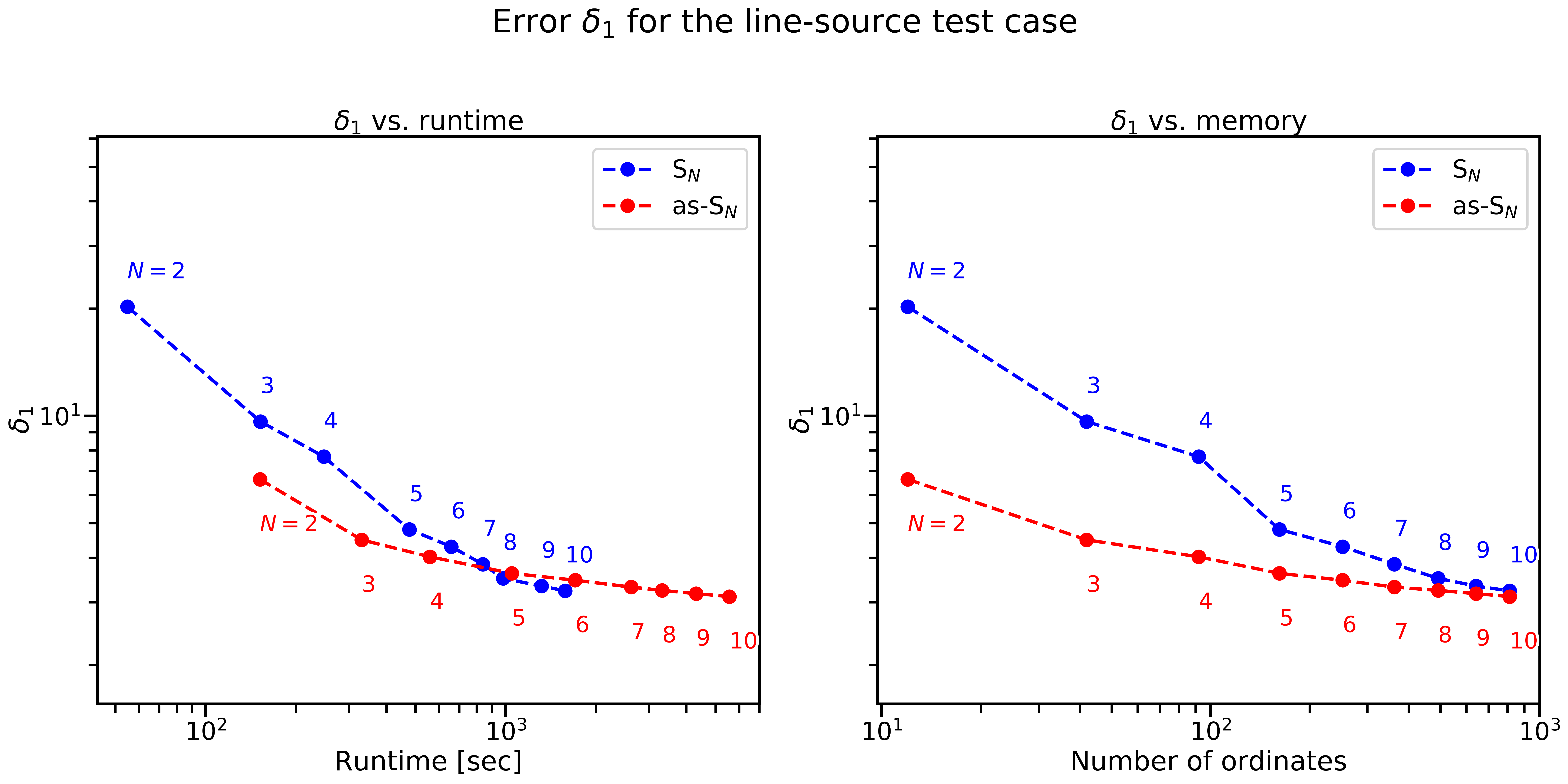}
	\caption{We computed $\delta_1 = \Vert\Phi_\text{numerical}-\Phi_\text{analytical}\Vert_2$ for the line-source test case using the implicit S$_N$ and as-S$_N$ method for $N_x \times N_y = 200\times 200$. Computations were performed on a quad core Intel\textsuperscript{\textregistered} i5-7300U CPU (2.60 GHz) with 12 GB memory.}
	\label{fig:rell2}
\end{figure}

\begin{figure}[h!]
	\centering
	\includegraphics[width=0.99\linewidth]{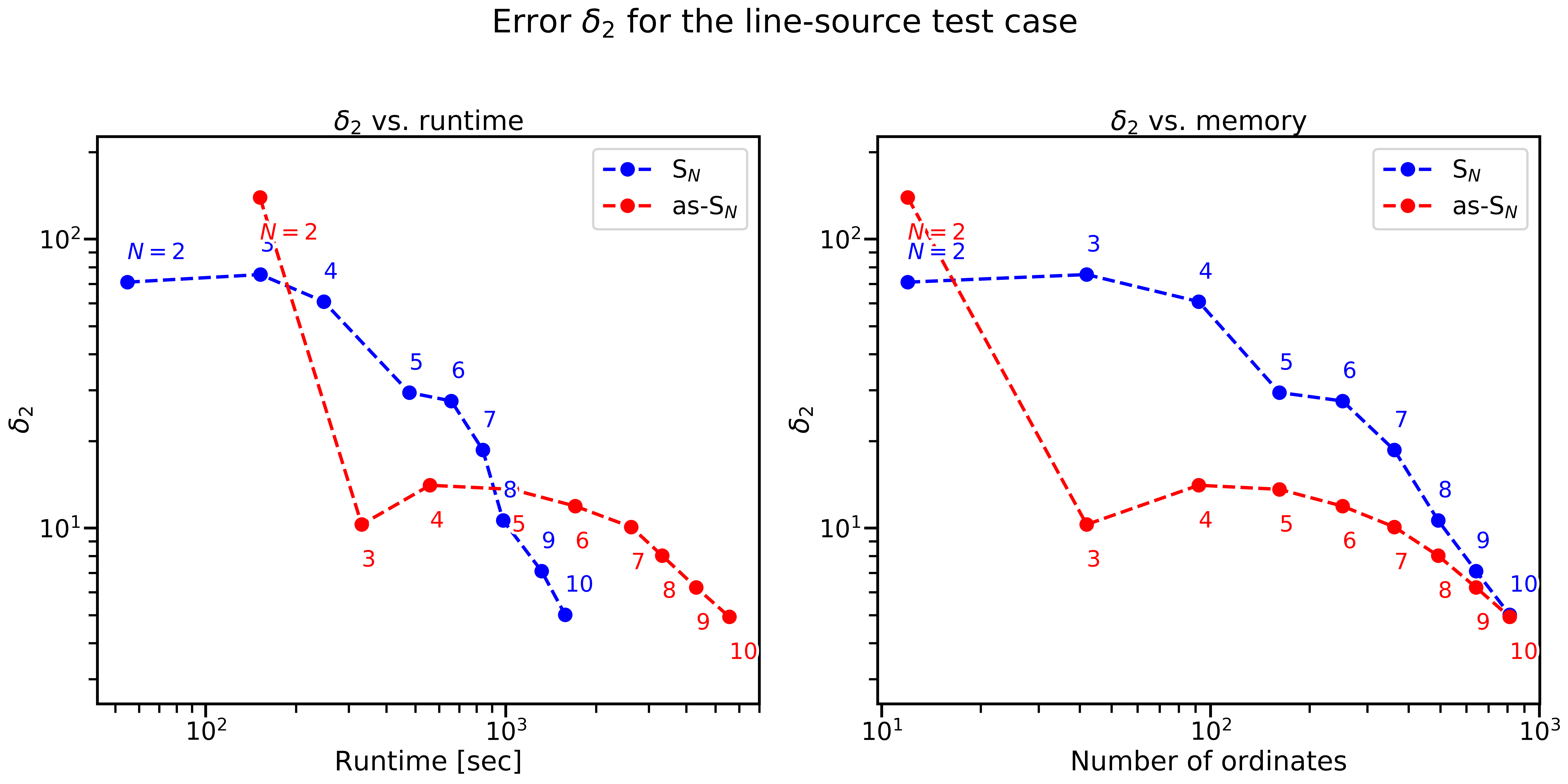}
\caption{We computed $\delta_2 = \Vert\nabla \Phi_\text{numerical}-\nabla\Phi_\text{analytical}\Vert_2$ for the line-source test case using the implicit S$_N$ and as-S$_N$ method for $N_x \times N_y = 200\times 200$. Computations were performed on a quad core Intel\textsuperscript{\textregistered} i5-7300U CPU (2.60 GHz) with 12 GB memory.}
	\label{fig:rell2grad}
\end{figure}

We also investigate the performance of the as-S$_N$ method when measured in runtime and in memory consumption. Consider therefore the results presented in Fig.~\ref{fig:rell2} and Fig.~\ref{fig:rell2grad}. Both figures summarize the results for the line-source test case computed with the \SN and \aSN methods for different values of $N$. Fig.~\ref{fig:rell2} measures the error between the numerical solution and the analytical solution in the L${}^2$ norm, called $\delta_1$. Fig.~\ref{fig:rell2grad} considers the $H^1$ semi-norm.

We observe an increase in runtime when activating artificial scattering, but a decrease in the errors $\delta_1$ and $\delta_2$. On the right, the errors are plotted against the number of ordinates which ultimately dictates the memory consumption. For example, an S$_8$ takes about as long as an as-S$_5$ computation and yields a similar $\delta_1$ error. However, the number of ordinates can be reduced from 492 to 162.
 For both, $\delta_1$ and $\delta_2$, the effect of artificial scattering vanishes in the limit of $N_q \rightarrow \infty$.

\subsection{Lattice test case}

We also investigate the lattice test case \cite{brunner2002forms,brunner2005two}, depicted in Fig.~\ref{fig:cb}. A constant, isotropic source is placed in the center of the domain in the orange square. In the white cells, the material is purely scattering, whereas the orange and black squares are purely absorbing. The boundary conditions are vacuum. All test case parameters are listed in Table \ref{tab:tabcheckerboard}.

\begin{minipage}{\textwidth}
	\begin{minipage}{0.35\textwidth}
		\includegraphics[width=0.99\linewidth]{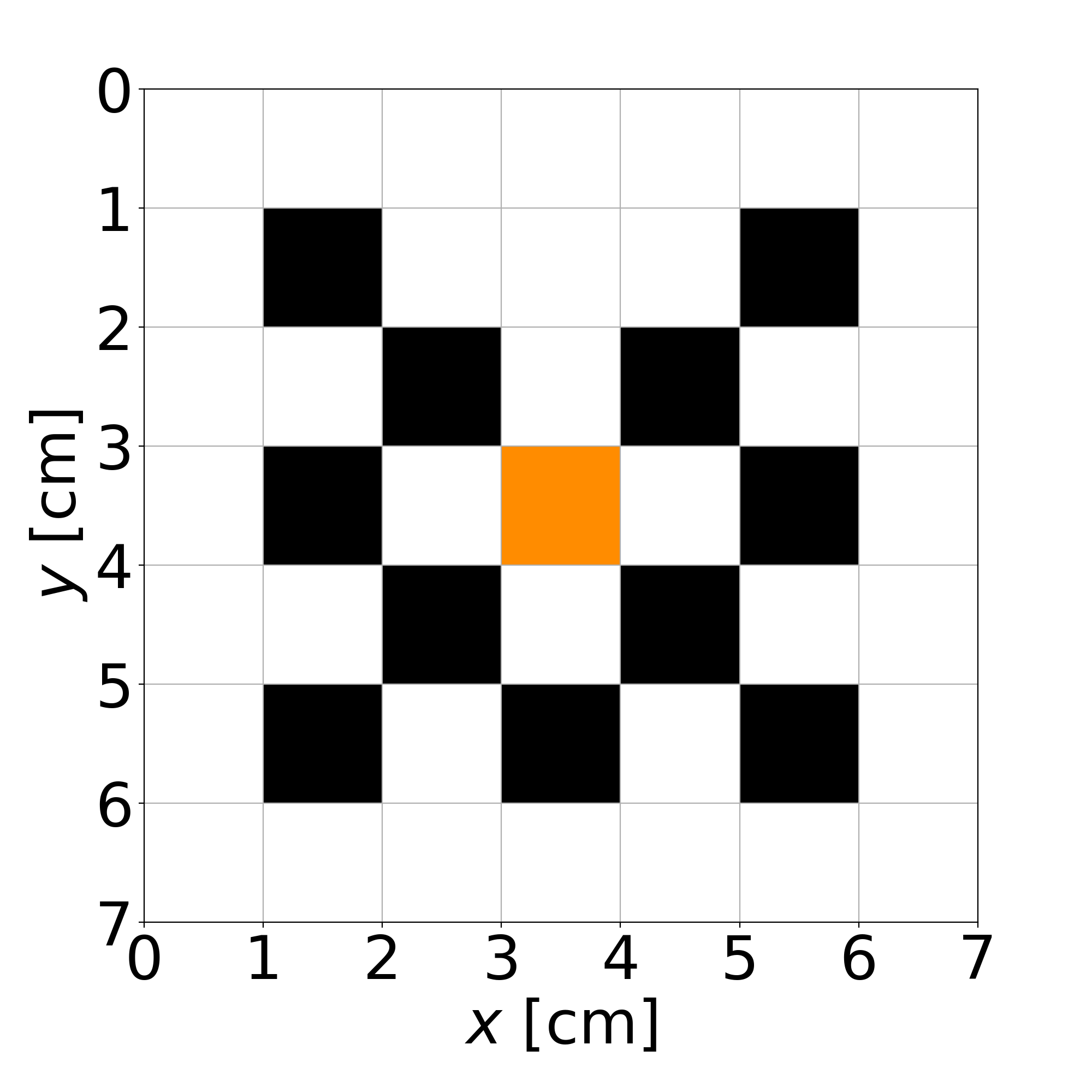}
		\captionof{figure}{Layout of the lattice test case.\\ Different materials (black, white, orange) \\with the source in the center (orange).}
		\label{fig:cb}
	\end{minipage}\hfill
	\begin{minipage}{0.55\textwidth}
		\begin{tabular}{c|c|c|c}
			Color & $\sigma_a$ in cm${}^{-1}$& $\sigma_s$ in cm${}^{-1}$ & $Q$ in cm${}^{-2}$s${}^{-1}$\\
			\hline
			\hline
			white & 0 & 1 & 0 \\
			black & 10 & 0 & 0 \\
			orange & 10 & 0 & 1 \\
		\end{tabular}
		\captionof{table}{Material properties for the lattice test case.\\ The domain is of size $[0\text{ cm},7\text{ cm}]^2$ and $t_{\text{end}}=3.2$ s. }
		\label{tab:tabcheckerboard}
	\end{minipage}
\end{minipage}

\begin{figure}[h!]
	\centering
	\includegraphics[width=0.9\linewidth]{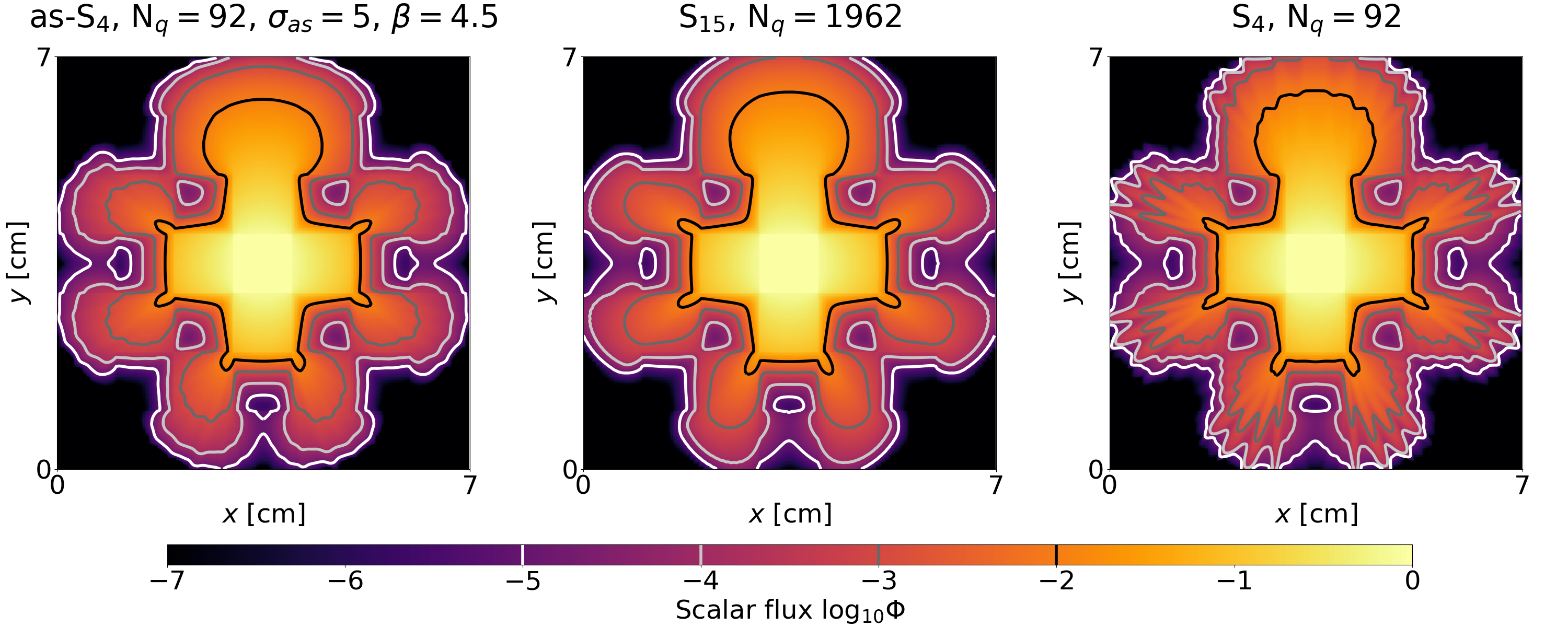}
	\caption{Comparison for the lattice test case at $t_{\text{end}}=3.2\, s$. Left the as-S$_4$ solution for the optimal parameter choice; center: the S$_{15}$ solution; right: the S$_{4}$ solution. Isolines are drawn at four different levels, highlighted inside the colorbar. We used $280\times280$ spatial cells.}
	\label{fig:checkerboard1}
\end{figure}

In Fig.~\ref{fig:checkerboard1} we see the as-S$_4$ solution to the lattice problem on the left, the S$_{15}$ solution in the center, and the S$_4$ solution on the right. Here, S$_{15}$ uses $1962$ ordinates while S$_4$ and as-S$_4$ use $92$ ordinates.

We take the S$_{15}$ solution with $N_q=1962$ as our reference solution. 
When comparing the as-S$_4$ solution with the S$_4$ solution, we see an improvement in the solution quality. Ray-effects are better mitigated in regions where the scalar flux is small. The number of ordinates is kept constant.

\begin{figure}[h!]
	\centering
	\includegraphics[width=0.9\linewidth]{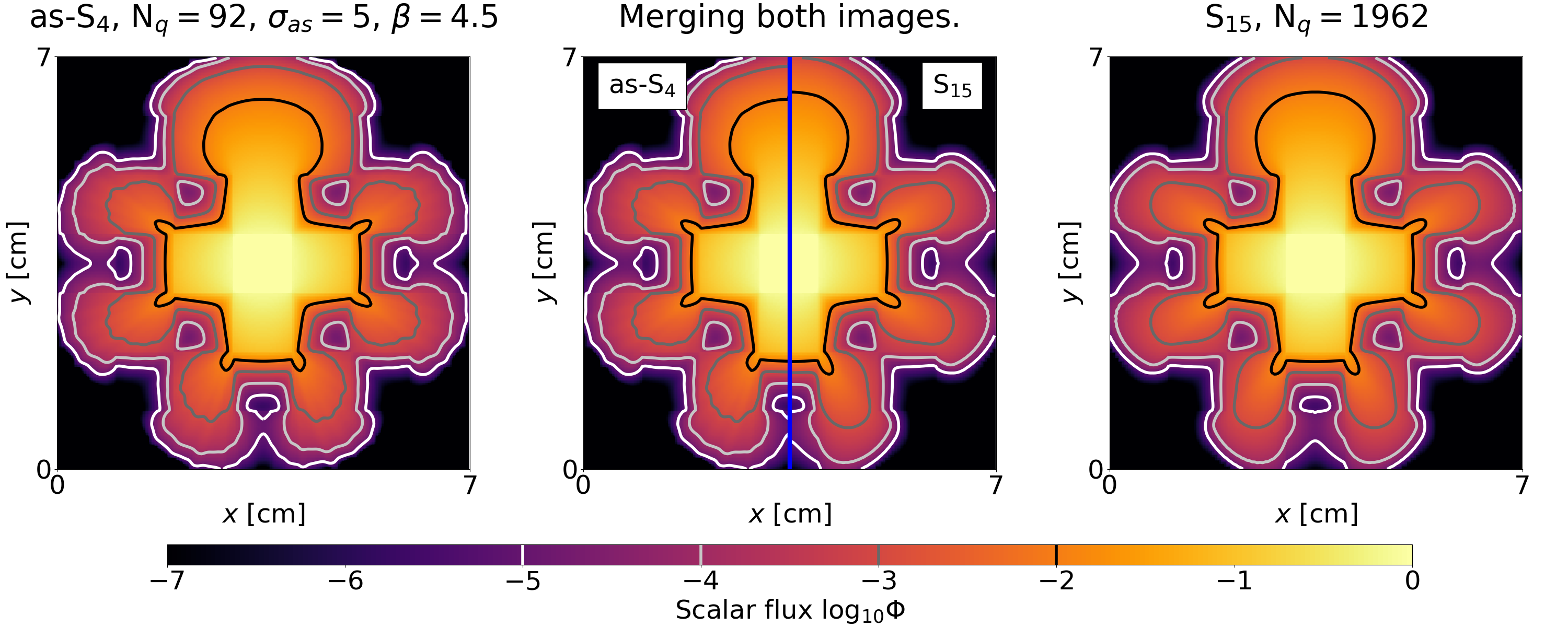}
	\caption{Comparison for the lattice test case at $t_{\text{end}}=3.2\, s$. Left: the as-S$_4$ solution; right: the S$_{15}$ solution; center: the image merges the left half of the left image with the right half of the right image. Isolines are drawn at four different levels, highlighted inside the colorbar. We used $280\times280$ spatial cells.}
	\label{fig:checkerboard2}
\end{figure}

Additionally, Fig.~\ref{fig:checkerboard2} puts the as-S$_4$ solution and the S$_{15}$ solution side-by-side in the center frame. 
Minor ray-effects are visible when looking at the white isoline.
However, the number of ordinates has been reduced by a factor of $\sim 21$.

Similar to the line-source test case, we set $\beta=4.5$ and $\sigma_{\text{as}}=5.0$ for explicit calculation, and $\beta=4.0$ and $\sigma_{\text{as}}=7.0$ for the implicit calculation.

\begin{figure}
	\centering
	\begin{subfigure}{.5\textwidth}
		\centering
		\includegraphics[width=0.99\linewidth]{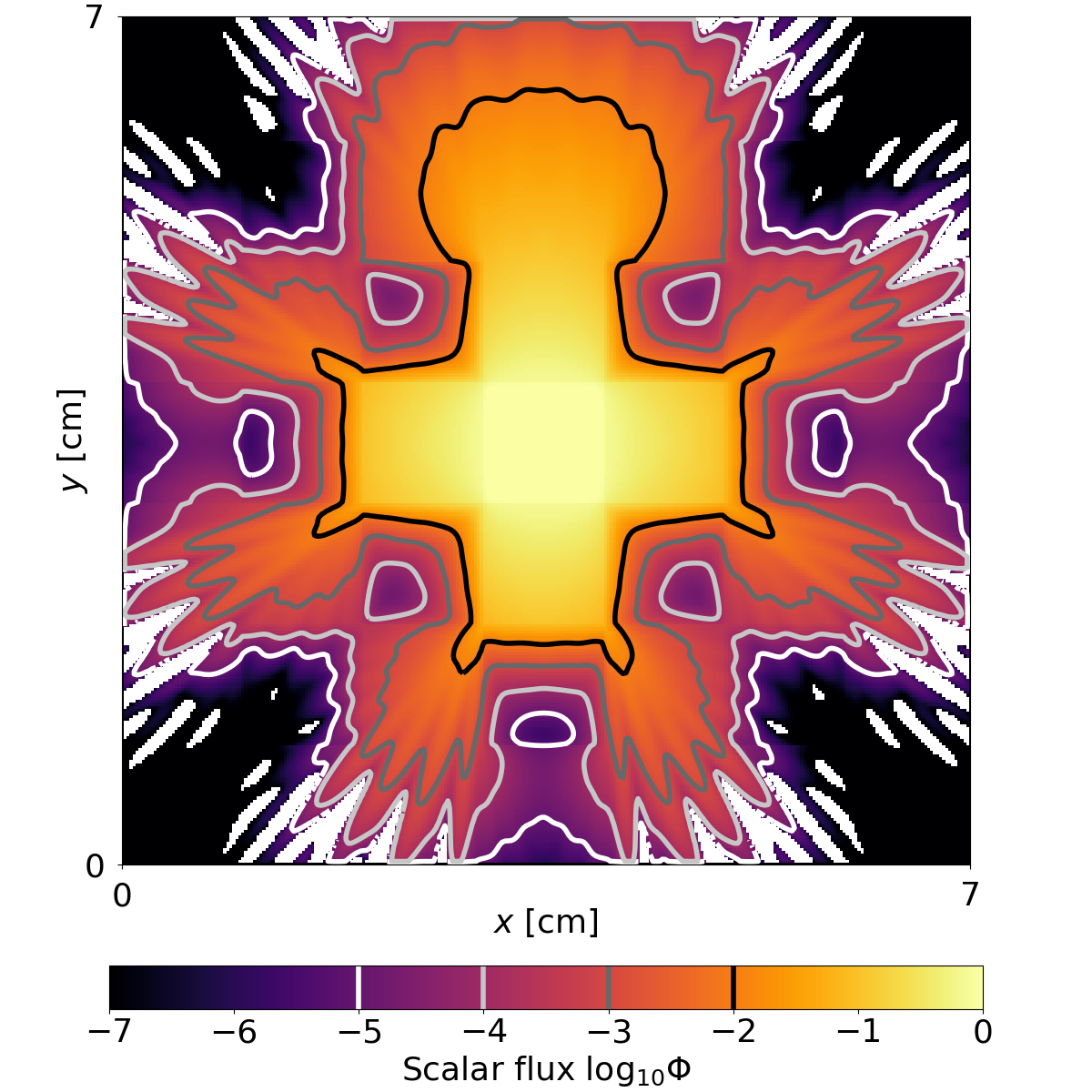}
		\caption{S$_4$.}
		\label{fig:cfl2}
	\end{subfigure}%
	\begin{subfigure}{.5\textwidth}
		\centering
		\includegraphics[width=0.99\linewidth]{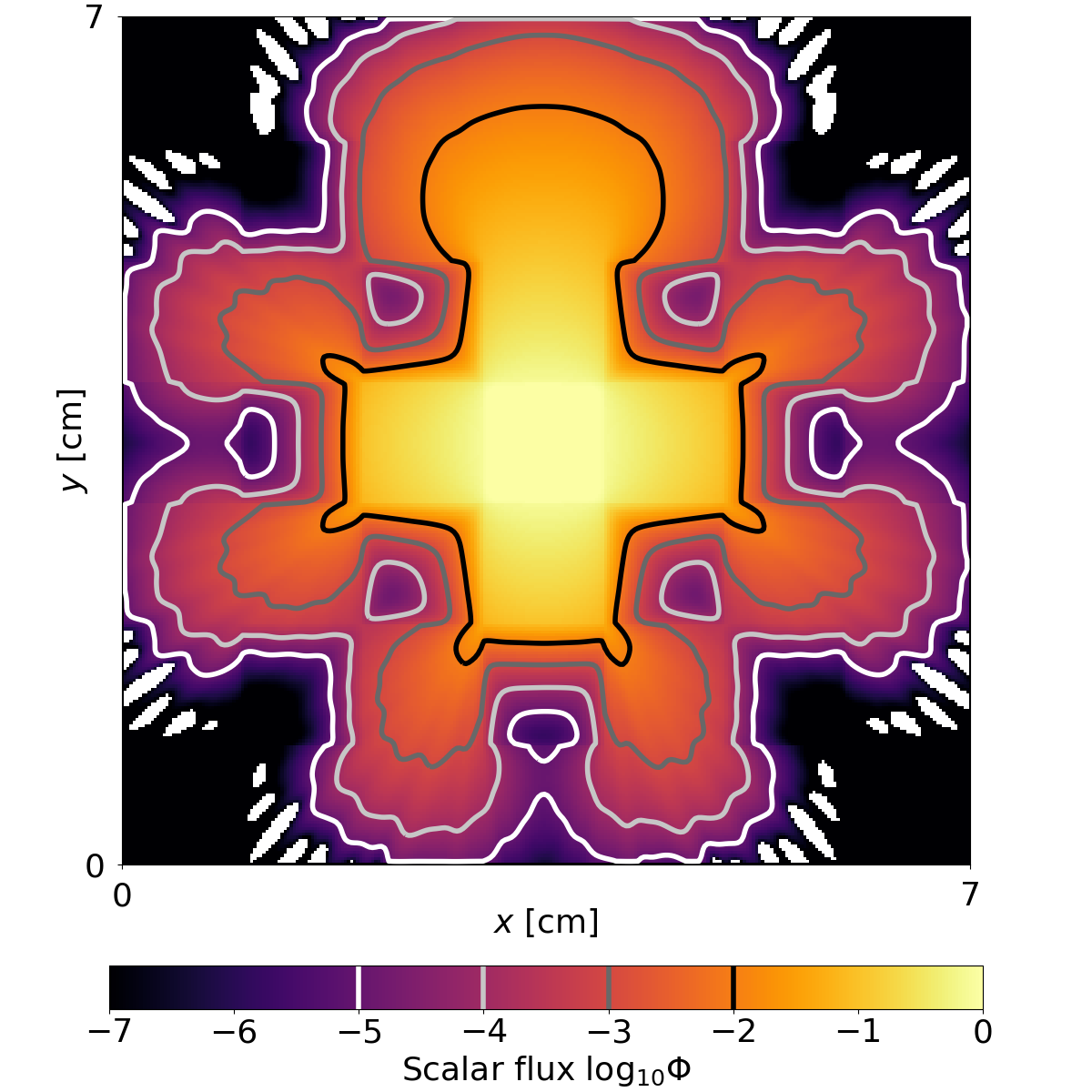}
		\caption{as-S$_4$ with $\sigma_{\text{as}}=7$ and $\beta=4$.}
		\label{fig:cfl2with}
	\end{subfigure}
	\caption{Solutions to the lattice problem with an implicit computation for a CFL number of $2$, $N_q=92$, and $N_x\times N_y=280\times 280$. The white regions indicate negativity of the solution.}
	\label{fig:lattice_impl_1}
\end{figure}

\begin{figure}%
	\centering
	\includegraphics[width=0.99\linewidth]{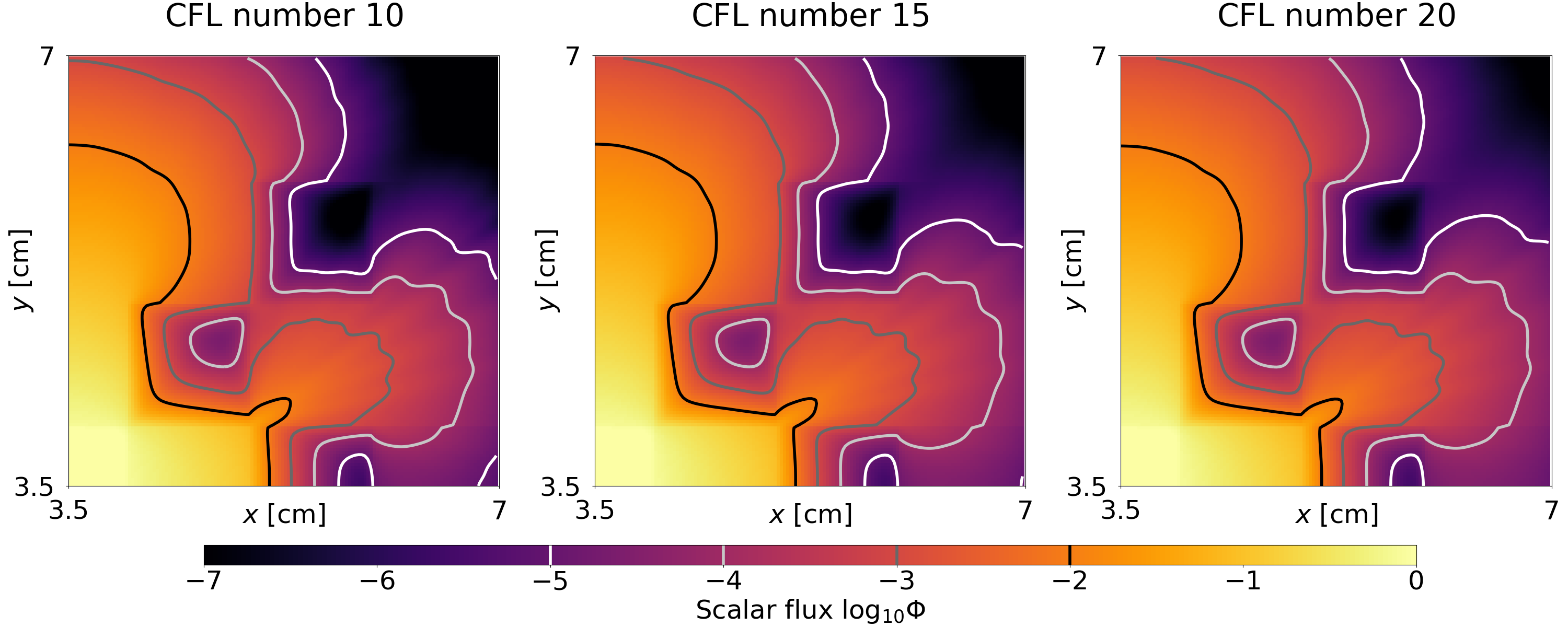}
	\caption{Solutions to the lattice problem with an implicit computation for different CFL numbers and $N_q=92$, $N_x\times N_y=280\times 280$, $\sigma_{\text{as}}=7$, and $\beta=4$. Zoom into the region  $[3.5,7]\times[3.5,7]$.}
	\label{fig:cfloverview}
\end{figure}

\FloatBarrier
 We also perform simulations for the lattice problem with the implicit time discretization. However, since the chosen scheme is only L${}^2$-stable, the solution becomes negative for the lattice test case as illustrated in Fig.~\ref{fig:lattice_impl_1}. Nevertheless, Fig.~\ref{fig:cfloverview} demonstrates the inherent advantage when performing implicit computations: We are able to use a very large CFL number, thus reducing the number of time steps and the overall computational costs drastically. Note that the scheme preserves positivity for the chosen CFL numbers.

\FloatBarrier
\section{Conclusion \& Outlook}
%
We have presented a new ray effect mitigation technique that relies on an additional, artificial scattering operator introduced into the radiative transfer equation. When the number of ordinates tends to infinity, the artificial scattering vanishes and the modified equation reduces to the original transport equation. In this case, when choosing the product of the scattering strength and the variance constant, the term tending to zero with the slowest rate is the Fokker-Planck operator.

The artificial scattering operator can be integrated into standard \SN codes. Solution algorithms, both for the explicit and implicit case, have been presented and rigorously analyzed in the non-standard, implicit case. To avoid using a large number of moments for the Krylov solver in the implicit case, we propose to invert the artificial scattering operator by a source iteration.

We have presented numerical results for the line-source and lattice test case. The results demonstrate that artificial scattering yields the same accuracy as \SN, but for a reduced number of ordinates.

For the second-order implicit computations the solutions might turn negative since L${}^2$ stability does not guarantee positivity of the solution. However, when choosing a sufficiently large CFL number, the solution values in our numerical experiments remain positive. A rigorous investigation of this effect and possibly the derivation of a CFL number ensuring positivity is left to future work.

Note that our test cases chose a constant value for the artificial scattering strength, however it seems plausible to make this strength spatially dependent to ensure that artificial scattering is only turned on when required. It remains to demonstrate the feasibility of the as-S$_N$ method in real-world applications using large-scale, highly parallelizable codes.

\section*{Acknowledgment}
The authors wish to thank Ryan G. McClarren (University of Notre Dame) for many fruitful discussions.
\setlength{\baselineskip}{12pt}
\bibliographystyle{siamplain}
\bibliography{main.bib}
\appendix
\section{Icosahedron quadrature}
The quadrature points and weights for the quadrature in Section \ref{subsec:quadrature} are given for order 2 (12 quadrature points). Every line contains four entries: the $x$, $y$, and $z$ position, as well as the quadrature weight. The quadrature weights sum to $4\pi$. All entries are  in \texttt{double} precision. The quadratures for order 2, order 3 (42 quadrature points), order 4 (92 quadrature points), and order 5 (162 quadrature points) can be downloaded as \texttt{.txt} files from a public repository at  \texttt{github.com/camminady/IcosahedronQuadrature}.

\begin{Verbatim}
+0.0000000000000000, +0.0000000000000000, +1.0000000000000000, +0.7370796188178727

+0.0000000000000000, -0.8944271909999159, +0.4472135954999580, +0.7370796188178727

+0.5257311121191336, -0.4472135954999580, +0.7236067977499790, +1.3293827143261787

-0.5257311121191336, -0.4472135954999580, +0.7236067977499790, +1.3293827143261787

+0.5257311121191336, +0.7236067977499790, +0.4472135954999580, +1.0751303204457423

-0.5257311121191336, +0.7236067977499790, +0.4472135954999580, +1.0751303204457423

+0.0000000000000000, +0.8944271909999159, -0.4472135954999580, +0.7370796188178727

-0.5257311121191336, +0.4472135954999580, -0.7236067977499790, +1.3293827143261787

+0.5257311121191336, +0.4472135954999580, -0.7236067977499790, +1.3293827143261787

+0.0000000000000000, +0.0000000000000000, -1.0000000000000000, +0.7370796188178727

-0.5257311121191336, -0.7236067977499790, -0.4472135954999580, +1.0751303204457423

+0.5257311121191336, -0.7236067977499790, -0.4472135954999580, +1.0751303204457423

\end{Verbatim}



\section{Implicit second order upwind scheme}
\label{app:upwind}
\begin{figure}[h!]
	\centering
	\includegraphics[width=0.95\linewidth]{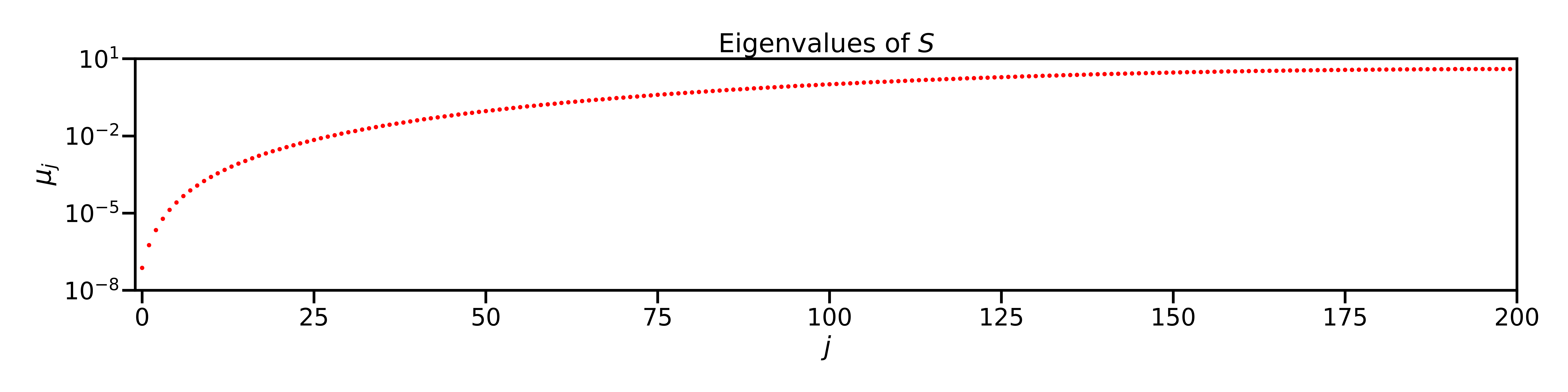}
	\caption{Eigenvalues of $S$ for $N_x=200$. All eigenvalues remain positive.}
	\label{fig:eigenvaluesS}
\end{figure}

In the following, we show that the chosen numerical flux is $L^2$ stable. For simplicity, we look at the one-dimensional advection equation
\begin{linenomath*}\begin{align}
\partial_t \psi + \Omega \partial_x \psi = 0
\end{align}\end{linenomath*}
with $\Omega\in\mathbb{R}_+$. A finite volume discretization is given by
\begin{linenomath*}\begin{align}\label{eq:scheme}
\psi_j^{n+1} = \psi_j^n - \lambda \left( g_{j+1/2} - g_{j-1/2}\right),
\end{align}\end{linenomath*}
where we use $\lambda := \Omega\Delta t / \Delta x$. A second order, implicit numerical flux is given by
\begin{linenomath*}\begin{align}
g_{j+1/2} := a \psi_j^{n+1} + b \psi_{j-1}^{n+1}
\end{align}\end{linenomath*}
with
\begin{linenomath*}\begin{align}
a := \frac32 , \enskip b:= -\frac12.
\end{align}\end{linenomath*}
Let us check if the scheme dissipates the $L^2$ entropy $\eta(\psi) := \psi^2/2$. For this we multiply our scheme \eqref{eq:scheme} with $\psi_j^{n+1}$, i.e. we obtain
\begin{linenomath*}\begin{align}\label{eq:term1}
\psi_j^{n+1}\psi_j^{n+1} = \psi_j^n \psi_j^{n+1} - \lambda \left( g_{j+1/2} - g_{j-1/2}\right)\psi_j^{n+1}.
\end{align}\end{linenomath*}
Now one needs to remove the cross term $\psi_j^n \psi_j^{n+1}$ which can be done by reversing the binomial formula
\begin{linenomath*}\begin{align}
\psi_j^n \psi_j^{n+1} = \frac12 (\psi_j^{n+1})^2 + \frac12 \left(\psi_j^{n}\right)^2 - \frac12 (\psi_j^{n+1}-\psi_j^{n})^2.
\end{align}\end{linenomath*}
Plugging this formulation for the cross term into \eqref{eq:term1} and making use of the definition of the square entropy $\eta$ gives
\begin{linenomath*}\begin{align}
\eta(\psi_j^{n+1}) = \eta(\psi_j^{n}) - \frac12 (\psi_j^{n+1}-\psi_j^{n})^2 - \lambda \left( g_{j+1/2} - g_{j-1/2}\right)\psi_j^{n+1}.
\end{align}\end{linenomath*}
This shows that in order to achieve entropy dissipation, i.e.
\begin{linenomath*}\begin{align}
\sum_{j=1}^{N_x}\eta(\psi_j^{n+1}) \leq \sum_{j=1}^{N_x}\eta(\psi_j^{n+1}),
\end{align}\end{linenomath*}
we need 
\begin{linenomath*}\begin{align*}\label{eq:dissipationL2Term}
\mathcal{E} = \sum_{j=1}^{N_x}\frac12 (\psi_j^{n+1}-\psi_j^{n})^2 + \lambda \sum_{j=1}^{N_x}\left( g_{j+1/2} - g_{j-1/2}\right)\psi_j^{n+1} \stackrel{!}{\geq} 0.
\end{align*}\end{linenomath*}
Note that the first term of $\mathcal{E}$, which essentially comes from the implicit time discretization is always positive. It remains to show that $\sum_{j=1}^{N_x}\left( g_{j+1/2} - g_{j-1/2}\right)\psi_j^{n+1}$ is positive as well. Let us rewrite this term for all spatial cells as a matrix vector product. I.e. when collecting the solution at time step $n+1$ for all $N_x$ spatial cells in a vector $\psi\in\mathbb{R}^{N_x}$, this term becomes 
\begin{linenomath*}\begin{align}
\sum_{j=1}^{N_x}\left(\left( g_{j+1/2} - g_{j-1/2}\right)\psi_j^{n+1}\right) = \psi^T B\psi
\end{align}\end{linenomath*}
where $B\in\mathbb{R}^{N_x\times N_x}$ is a lower triangular matrix. This product can be symmetrized with $S:=\frac{1}{2}(B+B^T)$, meaning that we have $\psi^T B\psi = \psi^T S\psi$. For our stencil, the matrix $S$ has entries $s_{jj} = \frac32$ on the diagonal and $s_{j,j-1} = s_{j-1,j} = -1$ as $s_{j,j-2} = s_{j-2,j} = \frac14$ on the lower and upper diagonals. Positivity of $\psi^T S\psi$ and thereby of the entropy dissipation term $\mathcal{E}$ in \eqref{eq:dissipationL2Term} is guaranteed if $S$ is positive definite, i.e. has positive eigenvalues. The eigenvalues for $S$ have been computed numerically to verify positivity in Fig.~\ref{fig:eigenvaluesS}.

\end{document}